\documentclass[11pt, reqno]{amsart}
\usepackage{amsfonts,latexsym}
\usepackage{amsmath}
\usepackage{amscd}
\usepackage{float,amsmath,amssymb,mathrsfs,bm,multirow,graphics}
\usepackage[dvips]{graphicx}
\usepackage[percent]{overpic}
\usepackage[numbers,sort&compress]{natbib}
\usepackage{enumerate}
\usepackage{amsaddr}

\newcommand{\nequation}{\setcounter{equation}{0}}

\newcommand{\R}{{\Bbb R}}

\newcommand{\C}{{\Bbb C}}

\newcommand{\Z}{{\Bbb Z}}
\newcommand{\proofbegin}{\noindent{\it Proof.\,\,}}
\newcommand{\proofend}{\hfill$\Box$\bigskip}
\newcommand{\proofendcontinue}{\hfill \raisebox{.8mm}[0cm][0cm]{$\bigtriangledown$}\bigskip}

\DeclareMathOperator{\tr}{tr}

\DeclareMathOperator{\im}{Im}
\DeclareMathOperator{\re}{Re}

\newcommand{\res}{\text{\upshape Res\,}}

\def\XXint#1#2#3{{\setbox0=\hbox{$#1{#2#3}{\int}$}
\vcenter{\hbox{$#2#3$}}\kern-.5\wd0}}

\allowdisplaybreaks

\newtheorem{theorem}{Theorem}
\newtheorem{proposition}{Proposition}[section]

\newtheorem{definition}[proposition]{Definition}

\newtheorem{remark}[proposition]{Remark}

\newtheorem{figuretext}{Figure}

\input epsf
\title[The defocusing nonlinear Schr\"odinger equation]{\sc The defocusing nonlinear Schr\"odinger equation with $t$-periodic data: New exact solutions}
    
\author{Jonatan Lenells}
\address{Department of Mathematics, KTH Royal Institute of Technology, \\ 100 44 Stockholm, Sweden.}
\email{jlenells@kth.se}

\begin{document}

\begin{abstract} 
\noindent
We consider solutions of the defocusing nonlinear Schr\"odinger (NLS) equation on the half-line whose Dirichlet and Neumann boundary values become periodic for sufficiently large $t$. We prove a theorem which, modulo certain assumptions, characterizes the pairs of periodic functions which can arise as Dirichlet and Neumann values for large $t$ in this way. The theorem also provides a constructive way of determining explicit solutions with the given periodic boundary values. Hence our approach leads to a class of new exact solutions of the defocusing NLS equation on the half-line. 
\end{abstract}

\maketitle

\noindent
{\small{\sc AMS Subject Classification (2010)}: 35Q55, 37K15.}

\noindent
{\small{\sc Keywords}: Initial-boundary value problem, time-periodic data, asymptotic behavior.}


\section{Introduction}\nequation
From the point of view of applications, one of the most important classes of initial-boundary value problems (IBVPs) consists of problems with asymptotically time-periodic boundary data. Such problems arise naturally, for example, when a wave entering a domain is measured at the domain's spatially fixed boundary.

Traditionally, nonlinear evolution PDEs have been studied mostly in the context of initial value problems. This has led to significant advances in our understanding of wellposedness issues and, at least for integrable equations, new solution generating techniques. Much less progress has been made in the context of IBVPs. For nonlinear integrable PDEs, the main obstacle in the analysis of IBVPs is that not all boundary values are known for a well-posed problem. Hence a succesful solution of the problem relies on the construction of the (generalized) Dirichlet to Neumann map. In general, the Dirichlet to Neumann map is highly nonlinear and can only be expressed in terms of a system of nonlinear integral equations \cite{BFS2003, F2005, trilogy1}. However, it has recently emerged that for asymptotically $t$-periodic data, the Dirichlet to Neumann map simplifies in the limit of large $t$, see \cite{BK2007, BIK2009, BKS2009, tperiodicI, tperiodicII}. In \cite{tperiodicII}, this observation was utilized to give an explicit construction of the Dirichlet to Neumann map for asymptotically $t$-periodic data in the limit of large $t$ and small data for the nonlinear Schr\"odinger (NLS) equation on the half-line. 

In this paper, we consider the defocusing NLS equation
\begin{align}\label{nls}
  iu_t + u_{xx} - 2 |u|^2 u = 0,
\end{align}
in the quarter plane  $\{x \geq 0, t \geq 0\}$ with asymptotically $t$-periodic data. The Dirichlet to Neumann map for (\ref{nls}) is a map $\{u_0(x), g_0(t)\} \mapsto g_1(t)$, where $u_0(x)$ denotes the initial data and $g_0(t) = u(0,t)$ and $g_1(t) = u_x(0,t)$ denote the Dirichlet data and the Neumann value, respectively. The time-periodicity implies that the Dirichlet to Neumann correspondence simplifies in the limit of large $t$ and it is crucial to take advantage of this simplification. Therefore, instead of considering the complete Dirichlet to Neumann map, we restrict attention to pairs of periodic functions $\{g_0^b(t), g_1^b(t)\}$ with the property that there exists a solution $u(x,t)$ on the half-line such that 
\begin{align}\label{uuxsim}
u(0,t) \sim g_0^b(t), \qquad u_x(0,t) \sim g_1^b(t), \qquad t \to \infty.
\end{align}
In many situations, the effect of the initial data diminishes in comparison to the periodic forcing as $t \to \infty$. Therefore a characterization of the periodic pairs $\{g_0^b(t), g_1^b(t)\}$ provides important information on the asymptotic behavior of the solution. 

Ideally, we would like to characterize all {\it asymptotically} admissible pairs, where a pair $\{g_0^b(t), g_1^b(t)\}$ is called asymptotically admissible if (\ref{uuxsim}) holds with polynomial convergence as $t \to \infty$.
Here we take a first step towards such a characterization by considering a characterization of the {\it eventually} admissible pairs, where a  pair $\{g_0^b(t), g_1^b(t)\}$ is called eventually admissible if there exists a $t_0 \geq 0$ and a solution $u(x,t)$ on the half-line such that
\begin{align}\label{ug0buxg1b}
u(0,t) = g_0^b(t), \qquad u_x(0,t) = g_1^b(t), \qquad t \geq t_0.
\end{align}
Clearly, every eventually admissible pair is asymptotically admissible.



Our main result (see Theorem \ref{mainth} below) characterizes, modulo certain assumptions, the periodic eventually admissible pairs for the defocusing NLS equation. Eventual admissibility is shown to be related to certain properties of the quotient $Q^b = B^b/A^b$, where the spectral functions $\{A^b(k), B^b(k)\}$ are defined in terms of the pair $\{g_0^b(t), g_1^b(t)\}$. 
Since the quotient $Q^b$ can be effectively computed from $\{g_0^b(t), g_1^b(t)\}$, the result provides a straightforward way of determining whether a given pair is eventually admissible. Moreover, the theorem provides a constructive way of finding an associated solution $u(x,t)$ with the given periodic boundary values. Hence the construction also leads to a class of new exact solutions of (\ref{nls}) on the half-line. The solutions in this class bear similarities with the stationary soliton solutions present for the focusing NLS. However, since the new solutions have singularities on the negative real axis, they are not regular solutions of (\ref{nls}) on the line, but only become regular solutions when restricted to the half-line $x \geq 0$. 

The proof of  Theorem \ref{mainth} utilizes the framework of \cite{tperiodicI}, which in turn relies on the ideas of the unified method introduced by Fokas \cite{F1997, F2002} as well as on the results for the focusing NLS with single exponential boundary values established by Boutet de Monvel and coauthors \cite{BK2007, BIK2009, BKS2009}. Single exponential pairs defined by
$$g_0^b(t) = \alpha e^{i\omega t}, \qquad g_1^b(t) = ce^{i\omega t}, \qquad t \geq 0,$$
where $\alpha>0$, $\omega \in \R$, and $c \in \C$, are of particular importance also for the defocusing NLS equation. As an application of our main result, we prove a theorem (see Theorem \ref{expth} below) which characterizes all eventually admissible single exponential pairs for the defocusing NLS equation. It turns out that $\{\alpha e^{i\omega t}, ce^{i\omega t}\}$ is eventually admissible if and only if $(\alpha, \omega, c)$ belongs to the set
$$\big\{(\alpha,\omega,c = - \alpha\sqrt{\omega + \alpha^2}) \; \big| \; \omega > 0, \; \alpha > 0\big\}.$$
In fact, by solving the associated RH problem explicitly, we find that if $(\alpha, \omega, c)$ belongs to this set, then 
\begin{align}\label{uexplicit}
u(x,t) = \frac{2 \alpha  \sqrt{\omega } \left(\sqrt{\alpha ^2+\omega }+\sqrt{\omega }\right) e^{x
   \sqrt{\omega }+i t \omega }}{\alpha ^2 \left(e^{2 x \sqrt{\omega }}-1\right)+2
   \sqrt{\omega } \left(\sqrt{\alpha ^2+\omega }+\sqrt{\omega }\right) e^{2 x
   \sqrt{\omega }}}, \qquad x \geq 0, \quad t \geq 0,
\end{align}
is an explicit solution of (\ref{nls}) in the quarter plane such that
$$u(0, t) = \alpha e^{i\omega t}, \qquad u_{x}(0,t) = c e^{i\omega t}, \qquad t\geq 0.$$
Note that although the solution (\ref{uexplicit}) possesses a singularity at
$$x = -\frac{\log \left(\frac{2 \sqrt{\omega } \sqrt{\alpha ^2+\omega }+\alpha ^2+2 \omega}{\alpha ^2}\right)}{2 \sqrt{\omega }} < 0,$$
it is smooth in the quarter plane $\{x \geq 0, t \geq 0\}$.

\section{Main results}\nequation\label{mainsec}

\begin{definition}\upshape\label{soldef}
A {\it solution of the NLS in the quarter plane} is a smooth function $u:[0,\infty) \times [0,\infty) \to \C$ such that $(a)$ $u(\cdot, t)$ belongs to the Schwartz class $\mathcal{S}([0,\infty))$ of rapidly decreasing functions for each $t \in [0, \infty)$, $(b)$ $u(x,t)$ satisfies (\ref{nls}) for $x> 0$ and $t > 0$, and $(c)$ $\|u(\cdot, t)\|_{L^1([0,\infty))}$ grows at most linearly in $t$ as $t \to \infty$.
\end{definition}

\begin{definition}\upshape\label{admissibledef}
A pair of smooth functions $\{g_0^b(t), g_1^b(t)\}$, $t \geq 0$, is {\it admissible} for NLS if there exists a solution $u(x,t)$ of the NLS in the quarter plane such that 
\begin{align}\label{g0g0Bg1g1B}
u(0, t) = g_0^b(t), \qquad u_{x}(0,t) = g_1^b(t), \qquad t\geq 0.
\end{align}
\end{definition}

\begin{definition}\upshape
A pair of smooth functions $\{g_0^b(t), g_1^b(t)\}$, $t \geq 0$, is {\it eventually admissible} for NLS if there exists a solution $u(x,t)$ of the NLS in the quarter plane and a $t_0 > 0$ such that 
\begin{align}\label{eventuallyadmissible}
u(0, t) = g_0^b(t), \qquad u_{x}(0,t) = g_1^b(t), \qquad t\geq t_0.
\end{align}
\end{definition}

Let $\{g_0^b(t), g_1^b(t)\}$ be a pair of smooth periodic functions of period $\tau = \frac{2\pi}{\omega} > 0$.
Following \cite{tperiodicI}, we define the entire $2 \times 2$-matrix valued function $Z(k)$ by $Z(k) = \psi(\tau,k)$ where $\psi(t,k)$ is the solution of the `background' $t$-part
\begin{align}\label{tpartB}
  \psi_t + 2ik^2\sigma_3 \psi = V^b\psi, \qquad \psi(0,k) = I,
\end{align}  
where
$$V^b(t,k) = \begin{pmatrix} -i |g_0^b(t)|^2 & 2kg_0^b(t) + ig_1^b(t) \\
2 k \bar{g}_0^b(t) - i \bar{g}_1^b(t) & i |g_0^b(t)|^2 \end{pmatrix}.$$
We define the spectral functions $A^b(k)$ and $B^b(k)$ by
\begin{align}\label{AbBbdef}
S^b(k) = \begin{pmatrix} \overline{A^b(\bar{k})} & B^b(k) \\ \overline{B^b(\bar{k})} & A^b(k) \end{pmatrix},
\end{align}
where
\begin{align*}
& S^b(k) = \sqrt{-\frac{Z_{11} - Z_{22} - \sqrt{G}}{2\sqrt{G}}} 
\begin{pmatrix} 1 & - \frac{2 Z_{12}}{Z_{11}-Z_{22}-\sqrt{G}} \\
 \frac{2 Z_{21}}{Z_{11}-Z_{22}-\sqrt{G}} & 1 \end{pmatrix}, 
 	\\
& G(k) = (\tr Z(k))^2 - 4.	
\end{align*} 
We are particularly interested in the quotient
\begin{align}\label{Qbdef}  
  Q^b(k) := \frac{B^b(k)}{A^b(k)} = - \frac{2 Z_{12}(k)}{Z_{11}(k) -Z_{22}(k) -\sqrt{G(k)}},
  \end{align}
and in the product
\begin{align}\label{Pbdef}  
  P^b(k) := \overline{B^b(\bar{k})} A^b(k) = - \frac{Z_{21}(k)}{\sqrt{G(k)}}.
  \end{align}
Individually, the functions $A^b$ and $B^b$ may have branch cuts generated by both of the square roots $\sqrt{-\frac{Z_{11} - Z_{22} - \sqrt{G}}{2\sqrt{G}}}$ and  $\sqrt{G}$. However, the combinations $Q^b$ and $P^b$ only involve $\sqrt{G}$, hence $Q^b$ and $P^b$ only have branch points at the zeros of $G(k)$ of odd order.

We can now state our main result.

\begin{theorem}\label{mainth}
  Let $\{g_0^b(t), g_1^b(t)\}$ be a pair of smooth periodic functions of period $\tau > 0$. 
  
  $(a)$ If $\{g_0^b(t), g_1^b(t)\}$ is eventually admissible for the defocusing NLS, then the function $Q^b(k) = \frac{B^b(k)}{A^b(k)}$ defined in (\ref{Qbdef}) satisfies the following properties:
  \begin{enumerate}[\upshape ({A}1)]
  \item $Q^b(k)$ is analytic for $\im k > 0$ and $Q^b(k)$ and all its derivatives have continuous extensions to $\im k \geq 0$.

  \item There exist complex constants $\{Q^b_j\}_1^\infty$ such that, for each $N \geq 1$,
  $$Q^b(k) = \frac{Q^b_1}{k} + \frac{Q^b_2}{k^2} +  \cdots + \frac{Q^b_{N-1}}{k^{N-1}} + O\Big(\frac{1}{k^N}\Big) \quad \text{uniformly as $k \to \infty$, $\im k \geq 0$.}$$
  Moreover, this expansion can be differentiated termwise to any order, i.e.
\begin{align}\label{dQbdkn}
\frac{d^nQ^b}{dk^n} = \frac{d^n}{dk^n}\bigg(\frac{Q^b_1}{k} + \frac{Q^b_2}{k^2} +  \cdots + \frac{Q^b_{N-1}}{k^{N-1}}\bigg) + O\Big(\frac{1}{k^{N}}\Big)
\end{align}
uniformly as $k \to \infty$, $\im k \geq 0$, for each $n \geq 1$ and each $N \geq 1$.

\item $\sup_{k \in \R} |Q^b(k)| < 1$.
\end{enumerate}

$(b)$ Suppose that $Q^b(k)$ satisfies {\upshape (A1)-(A3)} and that the (necessarily meromorphic) functions
\begin{align}\label{ZGquotients}
\text{$A^b(k)^2$ and $P^b(k)$ have at most finitely many poles.}
\end{align} 
Then, at least if the residues of $P^b$ are sufficiently small, the following statements hold:
\begin{itemize}
\item There exists a unique solution $u(x,t)$ of the defocusing NLS in the quarter plane which satisfies (\ref{g0g0Bg1g1B}); in particular, the pair $\{g_0^b(t), g_1^b(t)\}$ is admissible.

\item The solution $u(x,t)$ is given by
\begin{align}\label{ulim}
u(x,t) = 2i \lim_{k \to \infty} (kM(x,t,k))_{12}
\end{align}
where $M(x,t,k)$ is the unique solution of a RH problem (see (\ref{RHM}) below) formulated in terms of the function
$$h(k) = -\exp\bigg(\frac{1}{\pi i}\int_\R \frac{\log(1 - |Q^b(s)|^2)}{s-k} ds\bigg) P^b(k), \qquad \im k \geq 0.$$

\item The poles of $P^b(k)$ are simple and are contained in the set
\begin{align}\label{poleset}
 \Big\{\pm \frac{i\sqrt{n\omega}}{2} \, \Big| \, n = 1,2,\dots\Big\}.
\end{align}

\item The RH problem for $M(x,t,k)$ can be solved explicitly giving
\begin{align}\label{ulimMhat}
u(x,t) = 2i \lim_{k \to \infty} (k\hat{M}(x,t,k))_{12}
\end{align}
where $\hat{M}$ is defined as follows:
\begin{align}\nonumber
\hat{M}(x,t,k) = &\; (kI + B_N(x,t))(kI + B_{N-1}(x,t)) \times \cdots
	\\ \label{Mhatexplicit}
& \times (kI + B_1(x,t))\begin{pmatrix} \frac{1}{\prod_{j=1}^N(k-k_j)} & 0 \\ 0 & \frac{1}{\prod_{j=1}^N(k-\bar{k}_j)} \end{pmatrix},
\end{align}
where $\{k_j\}_1^N \subset i\R_+$ denote the simple poles of $P^b$ in the upper half plane and the matrices $\{B_j(x,t)\}_1^N$ are determined recursively from the algebraic equations
\begin{align}\label{Bjdef}
\begin{cases}
  (k_jI + B_j(x,t))\hat{M}_{j-1}(x,t,k_j) \begin{pmatrix} 1 \\ -d_j(x,t) \end{pmatrix} = 0, 
	\\
  (\bar{k}_jI + B_j(x,t))\hat{M}_{j-1}(x,t,\bar{k}_j) \begin{pmatrix} - \overline{d_j(x,t)} \\ 1 \end{pmatrix} = 0,
\end{cases} \qquad j = 1, \dots, N,
\end{align}
with 
\begin{align*}
& \hat{M}_0(x,t,k) = I,
	\\
& \hat{M}_j(x,t,k) = (kI + B_j(x,t))\hat{M}_{j-1}(x,t,k), \qquad j = 1, \dots, N-1,
	\\
& d_j(x,t) = -h_j \frac{\prod_{l=1, l \neq j}^N (k_j - k_l)}{\prod_{l=1}^N (k_j - \bar{k}_l)} e^{2i(k_jx + 2k_j^2t)}, \qquad j = 1, \dots, N,
	\\
& h_j = \underset{k_j}{\res} h(k), \qquad j = 1, \dots, N.
\end{align*}
\end{itemize}
\end{theorem}

      The proof of Theorem \ref{mainth} is presented in Section \ref{proof1sec}. 
      
      In the special case of single exponential profiles,  Theorem \ref{mainth} leads to the following result. 
      
\begin{theorem}[Single exponentials]\label{expth}
The pair $\{\alpha e^{i\omega t}, ce^{i\omega t}\}$ with $\alpha>0$, $\omega \in \R$, and $c \in \C$ is eventually admissible for the defocusing NLS if and only if $(\alpha, \omega, c)$ belongs to the set:
\begin{align}\label{reallyadmissibleD}
\big\{(\alpha,\omega,c = - \alpha\sqrt{\omega + \alpha^2}) \; \big| \; \omega > 0, \; \alpha > 0\big\}.
\end{align}
Moreover, if $(\alpha, \omega, c)$ belongs to (\ref{reallyadmissibleD}), then there is a unique solution $u(x,t)$ of the defocusing NLS in the quarter plane which satisfies 
$$u(0, t) = \alpha e^{i\omega t}, \qquad u_{x}(0,t) = c e^{i\omega t}, \qquad t\geq 0.$$
This solution is given explicitly by (\ref{uexplicit}).
\end{theorem}

      The proof of Theorem \ref{expth} is presented in Section \ref{proof2sec}.

\begin{remark}\upshape
The assumption in part $(b)$ of Theorem \ref{mainth} that the residues of $P^b$ are sufficiently small is only used to ensure that the RH problem (\ref{RHM}) has a unique solution for each $(x,t) \in [0,\infty) \times [0,\infty)$. If the solvability of the RH problem can be ensured in some other way (for example, by showing that the algebraic system (\ref{Bjdef}) is solvable for each $(x,t) \in [0,\infty) \times [0,\infty)$), then this assumption is not needed.
\end{remark}

\begin{remark}\upshape
  The condition that $\|u(\cdot, t)\|_{L^1([0,\infty))} = O(t)$ as $t \to \infty$ in Definition \ref{soldef} is only used once in the paper, namely, to ensure the validity of equation (\ref{GR}). 
\end{remark}

\section{Proof of Theorem \ref{mainth}}\label{proof1sec} \nequation

\subsection{Preliminaries}
Before we begin the proof, we collect some relevant definitions and facts. 
Given a pair of smooth periodic functions $\{g_0^b(t), g_1^b(t)\}$ of period $\tau > 0$, we define the entire function $Z(k)$ and $G(k)$ as in Section \ref{mainsec}. Then $Z(k)$ has eigenvalues $z(k)^{\pm1}$ where
$$z(k)^{\pm1} = \frac{1}{2} \big(\tr Z(k) \mp \sqrt{G(k)}\big).$$
We define the domains $\{D_j\}_1^4$ by
\begin{align}\nonumber
D_1 = \{\im k > 0\} \cap \{\im \tilde{\Omega}(k) > 0\},  \qquad
D_2 = \{\im k > 0\} \cap \{\im \tilde{\Omega}(k) < 0\}, 
	\\ \nonumber
D_3 = \{\im k < 0\} \cap \{\im \tilde{\Omega}(k) > 0\},  \qquad
D_4 = \{\im k < 0\} \cap \{\im \tilde{\Omega}(k) < 0\},
\end{align}
where $\tilde{\Omega}(k) = -\frac{\log z(k)}{i\tau}$.
We will also need the asymptotics of $Z(k)$ (see \cite{tperiodicI}):
\begin{align}\nonumber
& Z(k) = \begin{pmatrix} e^{-2ik^2\tau} & 0 \\ 0 & e^{2ik^2\tau} \end{pmatrix} 
+ \frac{1}{k}\begin{pmatrix} -i\eta_1(\tau)e^{-2ik^2\tau} & g_0^b(0)\sin(2k^2\tau) \\ 
 \overline{g_0^b(0)}\sin(2k^2\tau) & i\eta_1(\tau)e^{2ik^2\tau} \end{pmatrix}
	\\ \nonumber
& + \frac{1}{k^2} \begin{pmatrix} -i\overline{\eta_2(\tau)}e^{-2ik^2\tau}  -\frac{i}{2}|g_0^b(0)|^2 \sin(2k^2\tau) & \frac{ig_1^b(0)\sin(2k^2\tau)}{2} + g_0^b(0)\eta_1(\tau)\cos(2k^2\tau) \\ 
-\frac{i\overline{g_1^b(0)}\sin(2k^2\tau)}{2} + g_0^b(0)\eta_1(\tau)\cos(2k^2\tau) 
& i\eta_2(\tau)e^{2ik^2\tau} + \frac{i}{2}|g_0^b(0)|^2 \sin(2k^2\tau)
\end{pmatrix}
	\\\label{Zasymptotics}
& + O\bigg(\frac{e^{2ik^2\tau}}{k^3}\bigg) + O\bigg(\frac{e^{-2ik^2\tau}}{k^3}\bigg), \qquad k \to \infty, \quad k \in \C,
\end{align}
where $\eta_1(t) = \int_0^t \im(\overline{g_0^b(t')}g_1^b(t')) dt'$.
In particular,
\begin{align}\label{logzasymptotics}
\log z(k) = -2ik^2 \tau + O(k^{-1})
\end{align}
as $k \to \infty$ with $k$ remaining a bounded distance away from the branch cuts and the zeros of $\sin(2k^2\tau)$.

Given a contour $\Gamma \subset \C$ and $f \in L^2(\Gamma)$, we define the Cauchy transform $\mathcal{C}f$ by
\begin{align}\label{Cauchytransform}
(\mathcal{C} f)(k) = \frac{1}{2\pi i} \int_\Gamma \frac{f(s)}{s - k} ds, \qquad k \in \C \setminus \Gamma.
\end{align}
We let $\mathcal{C}_+ f, \mathcal{C}_- f \in L^2(\Gamma)$ denote the nontangential boundary values of $\mathcal{C} f$ as $k$ approaches $\Gamma$  from the left and right, respectively.

Suppose $\Gamma = \R$ and $f \in H^1(\R)$. Let $\Omega$ denote either the upper half-plane $\{\im k > 0\}$ or the lower half-plane $\{\im k < 0\}$. Then $\mathcal{C}f$ is holomorphic in $\Omega$ and extends to a uniformly H\"older continuous function of order $\frac{1}{2}$ on $\bar{\Omega}$, see Lemma 23.3 in \cite{BDT1988}; moreover, 
\begin{align}\label{supH1estimate}
\sup_{k \in \C \setminus \R} |(\mathcal{C}f)(k)| \leq \|f\|_{H^1(\R)} < \infty.
\end{align}

\subsection{Proof of $(a)$}
Suppose $\{g_0^b(t), g_1^b(t)\}$ is a periodic admissible pair and let $u(x,t)$ be an associated solution of the defocusing NLS in the quarter plane satisfying (\ref{g0g0Bg1g1B}).
As in \cite{tperiodicI}, we define the spectral functions $\{a(k),b(k), A(k),B(k)\}$ by
\begin{align*}
& s(k) = \begin{pmatrix}\overline{a(\bar{k})} & b(k) \\
   \overline{b(\bar{k})} & a(k) \end{pmatrix}, \qquad
  S(k) = \begin{pmatrix}\overline{A(\bar{k})} & B(k) \\
   \overline{B(\bar{k})} & A(k) \end{pmatrix},
\end{align*}
where  
$$s(k) = \mu_3(0,0,k), \qquad 	S(k) = \mu_1(0,0,k),$$
and $\mu_3(x, 0,k)$ and $\mu_1(0,t,k)$ are defined via linear Volterra integral equations in terms of $u(x,0)$ and $\{u(0,t), u_x(0,t), g_0^b(t), g_1^b(t)\}$, respectively.
Thanks to (\ref{g0g0Bg1g1B}), we have $S(k) = S^b(k)$, i.e.
$$A(k) = A^b(k), \qquad B(k) = B^b(k).$$
The spectral functions satisfy the global relation\footnote{The assumption that $\|u(\cdot, t)\|_{L^1([0,\infty))} = O(t)$ as $t \to \infty$ (see Definition \ref{soldef}) ensures that the global relation applies, see \cite{tperiodicI}.}
\begin{align}\label{GR}
 \frac{B^b(k)}{A^b(k)} = \frac{B(k)}{A(k)} = \frac{b(k)}{a(k)},
\end{align}
which is valid at least for all sufficiently large $k$ in the sector $\{\epsilon < \arg k < \frac{\pi}{2}- \epsilon\}$. Analytic continuation yields
$$Q^b(k) = \frac{b(k)}{a(k)}, \qquad \im k \geq 0.$$
Since $u(x,0)$ belongs to the Schwartz class $\mathcal{S}([0,\infty))$, the definition of $a(k)$ and $b(k)$ implies that $Q^b = b/a$ satisfies  (A1)-(A3), see \cite{FIS2005}. (Note that the function $a(k)$ is nonzero for $\im k \geq 0$ for the defocusing NLS.)
This proves $(a)$ of Theorem \ref{mainth} in the case that the pair $\{g_0^b(t), g_1^b(t)\}$ is admissible. The following result shows that the same implication holds also if the pair is eventually admissible.

\begin{proposition}
 Let $\{g_0^b(t), g_1^b(t)\}$ be smooth periodic functions of period $\tau > 0$. Then the pair $\{g_0^b(t), g_1^b(t)\}$ is eventually admissible for (\ref{nls}) if and only if it is admissible.
\end{proposition}
\proofbegin
Suppose $\{g_0^b(t), g_1^b(t)\}$ is an eventually admissible pair of $\tau$-periodic functions and let $u(x,t)$ be an associated solution of (\ref{nls}) in the quarter plane satisfying (\ref{eventuallyadmissible}). Choose $n \in \Z$ such that $n\tau > t_0$ and define a solutionÊ $v(x,t)$Ê of (\ref{nls}) by
$$v(x,t) = u(x,t+n\tau).$$
Then
$$v(0,t) = u(0,t+n\tau) = g_0^b(t+n\tau) = g_0^b(t), \qquad t \geq 0,$$
and
$$v_x(0,t) = u_x(0,t+n\tau) = g_1^b(t+n\tau) = g_1^b(t), \qquad t \geq 0.$$
The existence of the solution $v$ shows that the pair $\{g_0^b(t), g_1^b(t)\}$ is admissible.
The converse is trivial.
\proofend

\subsection{Proof of $(b)$}
Suppose $Q^b(k)$ satisfies (A1)-(A3).
Since $Z(k)$ and $G(k)$ are entire functions, property (A1) implies that $G$ has no zeros of odd order in $\im k \geq 0$. Indeed, such a zero together with the presence of the square root $\sqrt{G}$ would cause $Q^b$ not to be smooth in $\im k \geq 0$. 
Since $G(k) =\overline{G(\bar{k})}$, it follows that $G$ has no zeros of odd order in $\C$. Hence $\sqrt{G}$ is an entire function (this can be seen, for example, from the Hadamard factorization theorem). In particular, $(A^b)^2$, $Q^b$, and $P^b$ are meromorphic functions of $k \in \C$ and $z(k)$ is a nonzero entire function. It follows that $\tilde{\Omega}(k) = -\frac{\log z(k)}{i\tau}$ is an entire function. Liouville's theorem implies that $\tilde{\Omega}(k)$  is a polynomial of order $\leq 2$.
In fact, the asymptotic behavior (\ref{logzasymptotics}) of $\log z(k)$ shows that $\tilde{\Omega}(k) = 2k^2$. In particular, the domains $\{D_j\}_1^4$ coincide with the four quadrants of the complex plane. 
Moreover, since $z(k) = e^{-i\tau \tilde{\Omega}} = e^{-2ik^2\tau}$ we find
$$G(k) = -4\sin^2(2k^2\tau), \quad \text{i.e.} \quad \sqrt{G(k)} = 2i\sin(2k^2\tau).$$
It follows that
$$A^b(k)^2 = -\frac{Z_{11} - Z_{22} - \sqrt{G}}{2\sqrt{G}}, \qquad
P^b(k) = - \frac{Z_{21}(k)}{\sqrt{G(k)}},
$$
are meromorphic functions whose poles are contained in the set
$$\{k\,|\,\sin(2k^2\tau) = 0\} =  \Big\{\pm \frac{\sqrt{n\omega}}{2}, \pm \frac{i\sqrt{n\omega}}{2} \, \Big| \, n = 0, 1,2,\dots\Big\}.$$
All poles of $A^b(k)^2$ and $P^b(k)$ have order $\leq 1$ except for $k = 0$ which could be a double pole.

We henceforth suppose $A^b(k)^2$ and $P^b(k)$ have at most finitely many poles.

The remainder of the proof will proceed through a series of claims. 

\medskip\noindent
{\bf Claim 1.}  The function $a(k)$ defined by
\begin{align}\label{adef}
a(k) = \exp\bigg(-\frac{1}{2\pi i}\int_\R \frac{\log(1 - |Q^b(s)|^2)}{s-k} ds\bigg), \qquad \im k > 0,
\end{align}
has the following properties:
\begin{enumerate}[$(i)$]
  \item $a(k)$ is analytic for $\im k > 0$ and $a(k)$ and all its derivatives have continuous extensions to $\im k \geq 0$.

  \item There exist complex constants $\{a_j\}_1^\infty$ such that, for each $N \geq 1$,
\begin{align}\label{aexpansion} 
a(k) = 1 + \frac{a_1}{k} + \frac{a_2}{k^2} +  \cdots + \frac{a_{N-1}}{k^{N-1}} + O\Big(\frac{1}{k^{N}}\Big) \quad \text{uniformly as $k \to \infty$, $\im k \geq 0$.}
\end{align}
and this expansion can be differentiated termwise to any order.

  \item $a(k)$ is nonzero for $\im k \geq 0$.

  \item $|a(k)|^2 = \frac{1}{1 - |Q^b(k)|^2}$ for $k \in \R$.
\end{enumerate}

{\it Proof of Claim 1.}
Using the properties (A1)-(A3) of $Q^b$, property $(i)$ follows from (\ref{adef}) and the properties of the Cauchy operator $\mathcal{C}$ (see equation (\ref{supH1estimate})). Indeed, since $\log(1 - |Q^b|^2) \in H^1(\R)$, $a(k)$ is analytic for $\im k > 0$ with a continuous extension to $\im k \geq 0$. Differentiating (\ref{adef}), integrating by parts in the resulting expression, and recalling (\ref{dQbdkn}), we see that the same is true for $a'(k)$; the argument is easily extended to higher order derivatives of $a(k)$. 

We next prove $(ii)$. Let $N \geq 1$. By (A2), there exist real constants $\{r_j\}_2^N$ such that
$$\log(1 - |Q^b(s)|^2) = \frac{r_2}{s^2} + \frac{r_3}{s^3} + \cdots + \frac{r_N}{s^N} + O\bigg(\frac{1}{s^{N+1}}\bigg), \qquad |s| \to \infty, \quad s \in \R.$$
Let 
$$R_N(s) = \frac{r_2}{s^2} + \frac{r_3}{s^3} + \cdots + \frac{r_N}{s^N}.$$
Then
$$\int_{\R} \frac{R_N(s)}{s-k} ds = 0, \qquad \im k > 0.$$
Hence
$$a(k) = \exp\bigg(-\frac{1}{2\pi i}\int_\R \frac{\log(1 - |Q^b(s)|^2) - R_N(s)}{s-k} ds\bigg), \qquad \im k > 0.$$
The identity
$$\frac{1}{s-k} = -\sum_{j=1}^N \frac{s^{j-1}}{k^j} - \frac{s^N}{k^N(k - s)}$$
implies
\begin{align*}
a(k) =  \exp\bigg\{&\frac{1}{2\pi i} \sum_{j=1}^N \frac{1}{k^j}\int_\R s^{j-1}\big[\log(1 - |Q^b(s)|^2) - R_N(s)\big] ds
	\\
& + \frac{1}{2\pi i k^N} \int_\R s^N \frac{\log(1 - |Q^b(s)|^2) - R_N(s)}{k-s}  ds \bigg\}.
\end{align*}
In view of (\ref{supH1estimate}), this gives $(ii)$. Property $(iii)$ is obvious. 

Since $\mathcal{C}_+ - \mathcal{C}_- = I$, $\log a(k) = -\mathcal{C}\log(1 - |Q^b|^2)$ satisfies
$$\log a(k+i0) - \log a(k-i0) = -\log(1 - |Q^b(k)|^2), \qquad k \in \R.$$
But $\log \overline{a(\bar{k})} = - \log a(k)$ so $\log a(k-i0) = - \log \overline{a(k+i0)}$. This proves $(iv)$. 
\proofendcontinue

\medskip\noindent
{\bf Claim 2.} The function $b(k)$ defined by
\begin{align}\label{bdef}
b(k) = Q^b(k) a(k), \qquad \im k > 0,
\end{align}
has the following properties:
\begin{enumerate}[$(i)$]
  \item $b(k)$ is analytic for $\im k > 0$ and $b(k)$ and all its derivatives have continuous extensions to $\im k \geq 0$.

  \item There exist complex constants $\{b_j\}_1^\infty$ such that, for each $N \geq 1$,
\begin{align}\label{bexpansion} 
 b(k) = \frac{b_1}{k} + \frac{b_2}{k^2} +  \cdots + \frac{b_{N-1}}{k^{N-1}} + O\Big(\frac{1}{k^{N}}\Big) \quad \text{uniformly as $k \to \infty$, $\im k \geq 0$,}
\end{align}
and this expansion can be differentiated termwise to any order.

  \item $|a(k)|^2 - |b(k)|^2 = 1$ for $k \in \R$.
\end{enumerate}

{\it Proof of Claim 2.} 
Immediate from Claim 1 and the assumptions on $Q^b(k)$.
\proofendcontinue

\medskip\noindent
{\bf Claim 3.} The function $P^b(k)$ has the following properties:
\begin{enumerate}[$(i)$]
  \item There exist complex constants $\{P^b_j\}_1^\infty$ such that, for each $N \geq 1$,
\begin{align}\label{Pbexpansion} 
P^b(k) = \frac{P^b_1}{k} + \frac{P^b_2}{k^2} +  \cdots + \frac{P^b_{N-1}}{k^{N-1}} + O\Big(\frac{1}{k^{N}}\Big) \quad \text{uniformly as $k \to \infty$, $k \in \C$.}
\end{align}

  \item $P^b(k)$ is a rational function of $k \in \C$.  
\end{enumerate}

{\it Proof of Claim 3.} 
By assumption (\ref{ZGquotients}), $P^b$ has at most finitely many poles. If we can prove (\ref{Pbexpansion}), it will follow that $P^b(k)$ is a rational function.

By (\ref{Zasymptotics}),
\begin{align}\label{ZsqrtGasymptotics}
P^b(k) = -\frac{\frac{\overline{g_0^b(0)}\sin(2k^2\tau)}{k} + O\big(\frac{e^{2ik^2\tau}}{k^2}\big) + O\big(\frac{e^{-2ik^2\tau}}{k^2}\big)}{2i \sin(2k^2 \tau)}, \qquad k \to \infty, \quad k \in \C.
\end{align}
If $k \to \infty$ along a fixed ray $\arg k = \theta$ where $\theta$ is not a multiple of $\pi/2$, then (\ref{ZsqrtGasymptotics}) yields
$$P^b(k) = -\frac{\overline{g_0^b(0)}}{2ik} + O(k^{-2}), \qquad r \to \infty, \quad  k = re^{i\theta}.$$
Extending this argument to order $N$, we infer that (\ref{Pbexpansion}) holds as $k \to \infty$ along any fixed ray $\arg k = \theta$ disjoint from $\R \cup i\R$. 

The estimate (\ref{ZsqrtGasymptotics}) gives no information near the zeros $\{\pm \sqrt{\frac{n\pi}{2\tau}}, \pm i\sqrt{\frac{n\pi}{2\tau}} | n \geq 0\} \subset \R \cup i\R_+$ of $\sin(2k^2 \tau)$. In order to establish (\ref{Pbexpansion}) as $k \to \infty$ in a sector containing these zeros, we apply the Phragm\'en-Lindel\"of principle as follows.

Let $\{k_j\}$ denote the finite number of nonzero poles of $P^b = - Z_{21}/\sqrt{G}$ in $\C$. 
The functions $Z_{21}(k)k^2\prod_j (k-k_j)$ and $\sqrt{G(k)}$ are entire functions of order $\leq 2$, whose quotient is also entire; hence their quotient is also of order $\leq 2$ (see \cite{L1996}, p. 13). 

Choose $R > 0$ so large that $f(k) = k^N(P^b(k) - \sum_{n=1}^{N-1} \frac{P^b_n}{k^n})$ is holomorphic in the sector $S = \{|k| > R, -\pi/8 \leq \arg k \leq \pi/8\}$ and continuous on its closure $\bar{S}$. The function $f$ is bounded on $\partial S$ and there exists a $\rho < 4$ such that $|f(k)| \leq e^{C|k|^\rho}$ for $k \in S$. Consequently, by the Phragm\'en-Lindel\"of principle, $f(k)$ is bounded in $S$. This proves (\ref{Pbexpansion}) in $S$. The same argument applied to similar sectors $S$ which together cover all of $\C$ completes the proof of (\ref{Pbexpansion}). 
\proofendcontinue

\medskip\noindent
{\bf Claim 4.} The function $h(k)$ defined by
$$h(k) = - \frac{\overline{B^b(\bar{k})}A^b(k)}{a(k)^2} = - \frac{P^b(k)}{a(k)^2}, \qquad \im k \geq 0,$$
has the following properties:
\begin{enumerate}[$(i)$]
  \item $h(k)$ is meromorphic for $\im k > 0$ and $h(k)$ and all its derivatives extend continuously to $\im k \geq 0$. 
  
  \item The restriction of $h(k)$ to $\R$  is smooth and is given by
  \begin{align}\label{hba}
h(k) = -\frac{ \overline{b(k)}}{a(k)}, \qquad k \in \R.
\end{align}

  \item $h(k)$ has at most finitely many poles in the upper half-plane which all belong to the set 
 \begin{align}\label{hpoleset}
 \Big\{ \frac{i\sqrt{n\omega}}{2} \, \Big| \, n = 1,2,\dots\Big\}.
 \end{align}

  \item There exist complex constants $\{h_j\}_1^\infty$ such that, for each $N \geq 1$,
\begin{align}\label{hexpansion} 
h(k) = \frac{h_1}{k} + \frac{h_2}{k^2} +  \cdots + \frac{h_{N-1}}{k^{N-1}} + O\Big(\frac{1}{k^{N}}\Big) \quad \text{uniformly as $k \to \infty$, $\im k \geq 0$,}
\end{align}
and this expansion can be differentiated termwise to any order.

\end{enumerate}

{\it Proof of Claim 4.} 
We show (\ref{hba}); the remaining statements are then immediate from the properties of $P^b(k)$ and $a(k)$.

Since $\sqrt{G} = 2i\sin^2(2k^2\tau)$, $P^b = A^b \overline{B}^b$ and $Q^b = B^b/A^b$ are free from branch cuts. However, individually the functions $A^b(k)$  and $B^b(k)$  will still exhibit branch points at the possible zeros and poles of odd order of $A^b(k)^2$.
If there exists a branch cut that runs along the real axis, we denote by $A^b_+(k)$ and $A^b_-(k)$ the values of $A^b(k)$ on the top and bottom sides of the cut, respectively; similarly, $B^b_+(k)$ and $B^b_-(k)$ denote the values of $B^b(k)$ on the top and bottom sides of the cut. 
By (\ref{bdef}),
\begin{align}\label{GRquotient}
\frac{b(k)}{a(k)} = \frac{B^b(k)}{A^b(k)}, \qquad \im k > 0.
\end{align}
Taking the limit as $k$ approaches $\R$ from above and using that $Q^b$ has no branch cuts, we find
\begin{align}\label{baBpAp}
\frac{b(k)}{a(k)} = \frac{B^b_+(k)}{A^b_+(k)} = \frac{B^b_-(k)}{A^b_-(k)}, \qquad k \in \R.
\end{align}
Taking the same limit in the determinant relation $\det S^b = 1$, we obtain
$$A^b_+(k) \overline{A^b_-(k)} - B^b_+(k)\overline{B^b_-(k)} = 1, \qquad k \in \R.$$
Hence
\begin{align}\label{aaAA}
\frac{1}{|a(k)|^2} = 1 - \frac{b(k)}{a(k)} \frac{\overline{b(k)}}{\overline{a(k)}} 
= 1 - \frac{B^b_+(k)}{A^b_+(k)} \frac{\overline{B^b_-(k)}}{\overline{A^b_-(k)}}
= \frac{1}{A^b_+(k) \overline{A^b_-(k)}}, \qquad k \in \R.
\end{align}
This yields
$$h(k) 
= - \frac{(\overline{B^b(\bar{k})}A^b(k))_-}{a(k)^2}
= - \frac{\overline{B^b_+(k)}}{a(k)}\frac{A^b_-(k)}{a(k)}
= - \frac{\overline{B^b_+(k)}}{a(k)}\frac{\overline{a(k)}}{\overline{A^b_+(k)}}, \qquad k \in \R.$$
Together with (\ref{baBpAp}), this gives (\ref{hba}). 
\proofendcontinue

Equation (\ref{hba}) shows that $P^b(k) = \overline{b(k)}a(k)$ on $\R$. In particular, $P^b(k)$ has no poles in $\R$. This shows (\ref{poleset}).

\medskip\noindent
{\bf Claim 5.} The function $A^b(k)$ defined in (\ref{AbBbdef}) has the following properties:
\begin{enumerate}[$(i)$]  
  \item There exist complex constants $\{A^b_j\}_1^\infty$ such that, for each $N \geq 1$,
\begin{align}\label{Aexpansion} 
A^b(k)^2 = 1 + \frac{A^b_1}{k} + \frac{A^b_2}{k^2} +  \cdots + \frac{A^b_{N-1}}{k^{N-1}} + O\Big(\frac{1}{k^{N}}\Big) \quad \text{uniformly as $k \to \infty$, $k \in \C$.}
\end{align}

  \item $A^b(k)^2$ is a rational function of $k \in \C$.

\end{enumerate}

{\it Proof of Claim 5.} 
By assumption (\ref{ZGquotients}), $A^b(k)^2$ is a meromorphic function with at most finitely many poles. If we can prove (\ref{Aexpansion}), it will follow that $A^b(k)^2$ is a rational function.

By (\ref{Zasymptotics}),
\begin{align}\label{Z11Z22asymptotics}
A^b(k)^2 = -\frac{Z_{11} - Z_{22} - \sqrt{G}}{2\sqrt{G}} = \frac{2i\sin(2k^2\tau) + O\big(\frac{e^{2ik^2\tau}}{k}\big) + O\big(\frac{e^{-2ik^2\tau}}{k}\big)}{4i \sin(2k^2 \tau)}
+ \frac{1}{2}
\end{align}
as $k \to \infty$, $k \in \C$.
Thus, if $k \to \infty$ along a ray $\arg k = \theta$ where $\theta$ is not a multiple of $\pi/2$, then
$$A^b(k)^2 = 1 + O(k^{-1}), \qquad r \to \infty, \quad  k = re^{i\theta}.$$
Extending this argument to order $N$, we see that (\ref{Aexpansion}) holds as $k \to \infty$ along any such ray. The estimate (\ref{Z11Z22asymptotics}) gives no information near the zeros of $\sin(2k^2 \tau)$. 
However, as in the proof of (\ref{Pbexpansion}), we can appeal to the Phragm\'en-Lindel\"of principle to conclude that (\ref{Aexpansion}) holds.
\proofendcontinue

\begin{remark}\upshape
We have shown that $P^b = A^b \overline{B^b}$ and $(A^b)^2$ are rational functions. It follows that $Q^b = B^b/A^b$ also is a rational function. 
\end{remark}

\begin{figure}
\bigskip\bigskip
\begin{center}
\begin{overpic}[width=.45\textwidth]{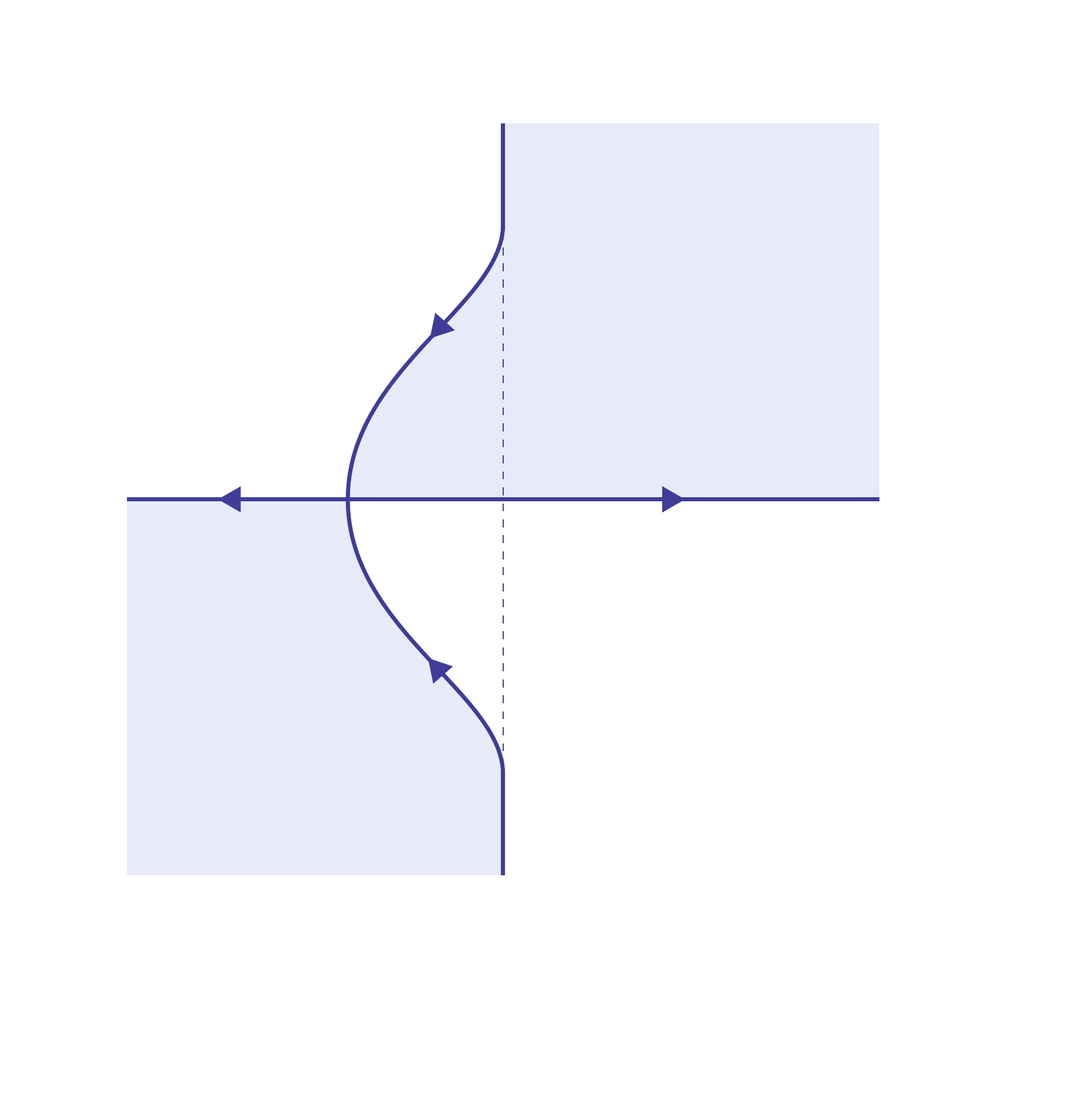}
      \put(70,75){$\mathcal{D}_1$}
      \put(20,75){$\mathcal{D}_2$}
      \put(20,20){$\mathcal{D}_3$}
      \put(70,20){$\mathcal{D}_4$}
      \put(70,42){$\Gamma$}
      \put(103,49){$\re k$}
      \put(44,103){$\im k$}
\end{overpic}
     \begin{figuretext}\label{Gammaoriented.pdf}
       The deformed quadrants $\{\mathcal{D}_j\}_1^4$ and the contour $\Gamma$.
      \end{figuretext}
     \end{center}
\end{figure}

Claims 4 and 5 imply that we may define deformed quadrants $\{\mathcal{D}_j\}_1^4$ in the complex $k$-plane such that no poles of $h(k)$, no poles or zeros of $A^b(k)^2$, and no branch cuts connecting such poles and zeros, lie in $\bar{\mathcal{D}}_2 \cup \bar{\mathcal{D}}_3$. 
We let $\Gamma = \R \cup (\bar{\mathcal{D}}_1 \cap \bar{\mathcal{D}}_2) \cup (\bar{\mathcal{D}}_3 \cap \bar{\mathcal{D}}_4)$ denote the contour separating the $\mathcal{D}_j$'s, oriented so that $\mathcal{D}_1$ and $\mathcal{D}_3$ lie to the left of $\Gamma$, see Figure \ref{Gammaoriented.pdf}. 
We choose the $\mathcal{D}_j$'s so that $\Gamma$ is invariant under the involution $k \to \bar{k}$. 

Consider the RH problem of finding a function $M(x,t,k)$ which $(i)$ is analytic for $k \in \C \setminus \Gamma$; $(ii)$ satisfies the asymptotic condition $M \to I$ as $k \to \infty$; and $(iii)$ satisfies the jump condition $M_+ = M_- J$ for $k \in \Gamma$ where
\begin{align}\label{Jdef}
&J(x,t,k) = \begin{cases} 
 \begin{pmatrix} 1 & 0 \\ h(k) e^{2i(kx + 2k^2t)} & 1 \end{pmatrix}, & k \in \bar{\mathcal{D}}_1 \cap \bar{\mathcal{D}}_2,
 	\\
I, & k \in \bar{\mathcal{D}}_2 \cap \bar{\mathcal{D}}_3,
	\\ 
 \begin{pmatrix} 1 & - \overline{h(\bar{k})} e^{-2i(kx + 2k^2t)} \\ 0 & 1 \end{pmatrix}, & k \in \bar{\mathcal{D}}_3 \cap \bar{\mathcal{D}}_4,
	\\
 \begin{pmatrix} 1 -  |h(k)|^2 & - \overline{h(\bar{k})} e^{-2i(kx + 2k^2t)} \\
h(k)e^{2i(kx + 2k^2t)} & 1 \end{pmatrix}, & k \in \bar{\mathcal{D}}_1 \cap \bar{\mathcal{D}}_4.
\end{cases}
\end{align}
We formulate this RH problem more precisely in the $L^2$-setting as follows:
\begin{align}\label{RHM}
\begin{cases}
M(x,t,\cdot) \in I + \dot{E}^2(\C\setminus \Gamma), 
	\\
M_+(x,t,k) = M_-(x, t, k) J(x, t, k) \quad \text{for a.e.} \ k \in \Gamma,
\end{cases}
\end{align}
where $\dot{E}^2(\C\setminus \Gamma)$ denotes the $L^2$-based generalized Smirnoff class; if $f$ is analytic in a domain $D$, then $f \in \dot{E}^2(D)$ provided there exist curves $\{C_n\}_1^\infty$ in $D$, tending to $\partial D$ in the sense that $C_n$ eventually surrounds each compact subset of $D$, such that $\sup_{n \geq 1} \|f\|_{L^2(C_n)} < \infty$.
The fact that $\det J = 1$ implies that if the RH problem (\ref{RHM}) has a solution, then it is necessarily unique. 

\begin{remark}\upshape
We refer to \cite{D1999, DZ2002a, DZ2002b, FIKN2006, LCarleson, Z1989} for more information on $L^2$-RH problems. Detailed definitions and proofs of all results used in this paper related to $L^2$-RH problems can be found in \cite{LCarleson}. 
\end{remark}

The rational function $P^b$ admits the partial fraction decomposition
$$P^b(k) = \sum_{i=1}^n \frac{r_i}{k - p_i}$$
where $p_i$, $i = 1, \dots, n$, denote the finitely many poles of $P^b(k)$ and $r_i$ are the corresponding residues.

\medskip\noindent
{\bf Claim 6.}  There exists an $\epsilon > 0$ such that if $|r_i| < \epsilon$ for each $i$, then the RH problem (\ref{RHM}) has a unique solution for each $(x,t) \in [0, \infty) \times [0,\infty)$.

{\it Proof of Claim 6.} 
By a small norm argument, the RH problem (\ref{RHM}) is uniquely solvable whenever $\|J - I\|_{L^\infty(\Gamma)}$ is small enough. Thus the unique solvability will follow if we can show that $h(k) = -P^b/a^2$ can be made arbitrarily small on $\Gamma$ by choosing the $r_i$ small. 

Since $(A^b)^2 ((A^b)^2 - 1) = (A^b)^2 B^b\overline{B^b} = P^b \overline{P^b}$ and $A^b \to 1$ as $k \to \infty$, we have 
$$A^b(k)^2 = \frac{1}{2}\Big(1 + \sqrt{1 + 4P^b(k) \overline{P^b(\bar{k})}}\Big), \qquad k \in \Gamma.$$ 
Thus by choosing the coefficients $r_i$ small, we can ascertain that the functions $P^b$ and $(A^b)^2 - 1$, and their derivatives, are arbitrarily small on $\Gamma$. But then $Q^b = P^b(A^b)^{-2}$ and its derivatives are also small on $\Gamma$. By (\ref{supH1estimate}) and (\ref{adef}) this implies that $a(k) - 1$ is small for $\im k \geq 0$. In particular, $a(k)$ is uniformly bounded away from zero. It follows that by choosing the $r_i$'s small, we can ascertain that $h(k)$ is arbitrarily small on $\Gamma$.
\proofendcontinue

In the remainder of the proof we assume that the RH problem (\ref{RHM}) has a solution for each $(x,t) \in [0, \infty) \times [0,\infty)$ (as a consequence of the residues of $P^b$ being small or for some other reason).

\medskip\noindent
{\bf Claim 7.}  The limit
\begin{align}\label{ulim2}
u(x,t) = 2i \lim_{k \to \infty} (kM(x,t,k))_{12}
\end{align}
exists for every $(x,t) \in [0,\infty) \times [0,\infty)$, and the function $u(x,t)$ defined by (\ref{ulim2}) is a smooth function of $(x,t) \in [0,\infty) \times [0,\infty)$ which satisfies the defocusing NLS equation for $x \geq 0$ and $t \geq 0$.

\begin{figure}
\bigskip\bigskip
\begin{center}
\begin{overpic}[width=.45\textwidth]{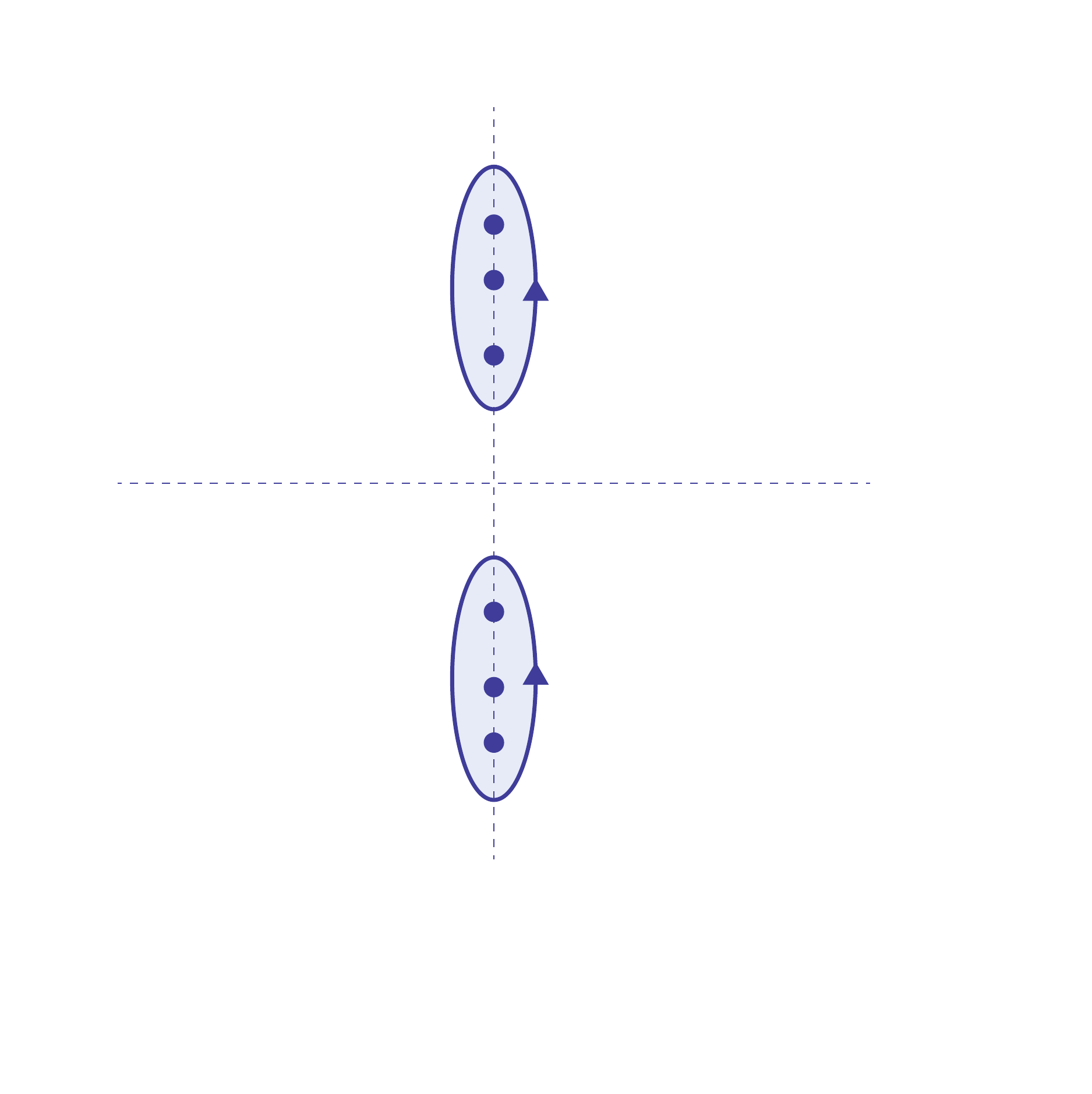}
      \put(58,76){$\gamma_1$}
      \put(58,23){$\gamma_2$}
      \put(102.5,48.5){$\re k$}
      \put(44,103.5){$\im k$}
\end{overpic}
     \begin{figuretext}\label{gamma.pdf}
       The contour $\gamma = \gamma_1 \cup \gamma_2$ in the complex $k$-plane. The dots indicate the poles of $h(k)$ and $\overline{h(\bar{k})}$. The open set $U$ is the shaded area enclosed by $\gamma$. 
      \end{figuretext}
     \end{center}
\end{figure}

{\it Proof of Claim 7.}
Let $\gamma_1$ be a contour in $\mathcal{D}_1$ that surrounds the poles of $h(k)$ in the upper half-plane in the counterclockwise direction, see Figure \ref{gamma.pdf}. Let $\gamma_2$ be the complex conjugate of $\gamma_1$ also oriented in the counterclockwise direction. Let $\gamma = \gamma_1 \cup \gamma_2$ and let $U$ denote the open region enclosed by $\gamma$. 
Then $h \in (\dot{E}^2 \cap E^\infty)(\mathcal{D}_1 \setminus \bar{U})$, where $E^\infty(D)$ denotes the set of bounded analytic functions in a domain $D$. Hence, by a contour deformation argument (see e.g. Section 5 of \cite{LCarleson}), the function $m(x,t,k)$ defined by
$$m = \begin{cases} 
M \begin{pmatrix} 1 & 0 \\ -h(k) e^{2i(kx + 2k^2t)} & 1 \end{pmatrix}, & k \in \mathcal{D}_1 \setminus \bar{U},
	\\
M \begin{pmatrix} 1 & -\overline{h(\bar{k})} e^{-2i(kx + 2k^2t)} \\ 0 & 1 \end{pmatrix}, & k \in \mathcal{D}_4 \setminus \bar{U},
	\\
M, &	\text{otherwise},
\end{cases}$$
satisfies the $L^2$-RH problem
\begin{align}\label{RHm}
\begin{cases}
m(x,t,\cdot) \in I + \dot{E}^2(\C \setminus \gamma), 
	\\
m_+(x,t,k) = m_-(x, t, k) v(x, t, k) \quad \text{for a.e.} \ k \in \gamma,
\end{cases}
\end{align}
where
\begin{align*}
&v(x,t,k) = \begin{cases} 
 \begin{pmatrix} 1 & 0 \\ h(k) e^{2i(kx + 2k^2t)} & 1 \end{pmatrix}, & k \in \gamma_1,
	\\ 
 \begin{pmatrix} 1 & \overline{h(\bar{k})} e^{-2i(kx + 2k^2t)} \\ 0 & 1 \end{pmatrix}, & k \in \gamma_2.
\end{cases}
\end{align*}
It follows that $w(x,t,\cdot) := v(x,t,\cdot) - I$ is a smooth function of $(x,t) \in [0,\infty)\times [0,\infty)$ into $L^p(\gamma)$ for each $1 \leq p \leq \infty$. We define the operator $\mathcal{C}_w: L^2(\gamma) + L^\infty(\gamma) \to L^2(\gamma)$ by 
\begin{align}\label{Cwdef}
\mathcal{C}_w(f) = \mathcal{C}_-(f w).
\end{align}
Since $w$ is nilpotent, our assumption that the RH problem has a unique solution implies that $I - \mathcal{C}_w$ is a bijective map $L^2(\gamma) \to L^2(\gamma)$ and that
$$m(x,t,k) = I + \frac{1}{2\pi i} \int_{\gamma} \frac{(\mu w)(x,t, s)ds}{s-k}, \qquad k \in \C \setminus \gamma,$$
where $\mu = I + (I - \mathcal{C}_w)^{-1}\mathcal{C}_wI \in I + L^2(\gamma)$ cf. \cite{Z1989} (see Section 5 of \cite{LCarleson} for a detailed proof). 
Hence the limit 
\begin{align}\nonumber
u(x,t)& = 2i \lim_{k \to \infty} (kM(x,t,k))_{12} 
= 2i \lim_{k \to \infty} (km(x,t,k))_{12}
	\\ \label{ulim3}
 &= -\frac{1}{\pi}\int_\gamma (\mu w)_{12}(x,t,s) ds
\end{align}
exists for each $(x,t)$. Since $(x,t) \mapsto w(x,t,\cdot)$ is smooth and $I - \mathcal{C}_w$ is invertible for each $(x,t)$, $(x,t)\mapsto \mu(x,t,\cdot) - I$ is a smooth map $[0,\infty) \times [0,\infty) \to L^2(\gamma)$. It follows from (\ref{ulim3}) that $u(x,t)$ is a smooth function of $(x,t) \in [0, \infty) \times [0,\infty)$. 

The fact that $u$ satisfies (\ref{nls}) follows from the dressing method. 
\proofendcontinue

The deformed quadrants $\mathcal{D}_j$ are defined so that
\begin{align}\label{hApm1}
h \in \dot{E}^2(\mathcal{D}_2) \cap E^\infty(\mathcal{D}_2), \qquad
A^b,\frac{1}{A^b} \in 1 + \dot{E}^2(\mathcal{D}_2\cup \mathcal{D}_3) \cap E^\infty(\mathcal{D}_2 \cup \mathcal{D}_3).
\end{align}

\medskip\noindent
{\bf Claim 8.}  $u(\cdot,0) \in \mathcal{S}([0,\infty))$.

{\it Proof of Claim 8.}
Define $u_0(x)$ by
\begin{align}\label{u0lim}
u_0(x) = 2i \lim_{k\to\infty} (km^{(x)}(x,k))_{12},
\end{align}
where $m^{(x)}$ satisfies the $L^2$-RH problem
\begin{align}\label{mxRH}
\begin{cases}
m^{(x)}(x, \cdot) \in I + \dot{E}^2(\C \setminus \R),\\
m_+^{(x)}(x,k) = m_-^{(x)}(x, k) \begin{pmatrix} \frac{1}{|a(k)|^2} &  \frac{b(k)}{\overline{a(\bar{k})}} e^{-2ikx} \\
-\frac{\overline{b(\bar{k})}}{a(k)} e^{2ikx} & 1 \end{pmatrix},  \qquad k \in \R.
\end{cases}
\end{align}
The unique solvability of this RH problem follows from the existence of a vanishing lemma \cite{Z1989}. The smoothness and decay properties of $\{a(k),b(k)\}$ derived in Claims 1 and 2 imply that $u_0(x)$ is well defined by (\ref{u0lim}) and belongs to the Schwartz class $\mathcal{S}([0,\infty))$ cf. \cite{FIS2005}. 

On the other hand, by (\ref{hApm1}),
$$h(k) e^{2ikx} \in \dot{E}^2(\mathcal{D}_2) \cap E^\infty(\mathcal{D}_2), \qquad \overline{h(\bar{k})} e^{-2ikx} \in \dot{E}^2(\mathcal{D}_3) \cap E^\infty(\mathcal{D}_3).$$ 
Thus a contour deformation argument together with (\ref{aaAA}) and the expression (\ref{Jdef}) for $J$ imply that the function $M^{(x)}(x,k)$ defined by
$$M^{(x)}(x,k) = \begin{cases} M(x, 0, k), & k \in \mathcal{D}_1 \cup \mathcal{D}_4,
	\\
M(x,0,k) \begin{pmatrix} 1 & 0 \\ h(k) e^{2ikx} & 1 \end{pmatrix}, & k \in \mathcal{D}_2,
	\\
M(x,0,k) \begin{pmatrix} 1 & \overline{h(\bar{k})} e^{-2ikx} \\ 0 & 1 \end{pmatrix}, & k \in \mathcal{D}_3,
\end{cases}$$
also satisfies (\ref{mxRH}). By uniqueness, $m^{(x)} = M^{(x)}$. Comparing the definition (\ref{ulim}) of $u(x,t)$ with (\ref{u0lim}), we obtain $u(x,0) = u_0(x)$ for $x \geq 0$.
\proofendcontinue

\medskip\noindent
{\bf Claim 9.}  $u(0,t) = g_0^b(t)$ and $u_x(0,t) = g_1^b(t)$ for $t \geq 0$.

{\it Proof of Claim 9.}

Define $m^{(t)}(t,k)$ in terms of $g_0^b$ and $g_1^b$ by
$$m^{(t)}(t,k) = \begin{cases}  ([\psi(t,k)]_1e^{2ik^2t}, \frac{[\mathcal{E}(t,k)]_2}{A^b(k)}), \qquad k \in \mathcal{D}_1 \cup \mathcal{D}_3, \\ (\frac{[\mathcal{E}(t,k)]_1}{\overline{A^b(\bar{k})}}, [\psi(t,k)]_2e^{-2ik^2t}), \qquad k \in \mathcal{D}_2 \cup \mathcal{D}_4. \end{cases}$$
where $\psi(t,k)$ is the solution in (\ref{tpartB}) and
$$\mathcal{E}(t,k) = \psi(t,k) S^b(k) e^{2ik^2t\sigma_3}.$$
Since
$$\frac{[\mathcal{E}(t,k)]_1}{\overline{A^b(\bar{k})}} = \big([\psi(t,k)]_1 + [\psi(t,k)]_2 \overline{Q^b(\bar{k})}\big)e^{2ik^2t},$$
property (A1) shows that $m^{(t)}(t, \cdot)$ is analytic in $\mathcal{D}_4$. 
Equation (\ref{hApm1}) implies that $m^{(t)}(t, \cdot)$ is analytic in $\mathcal{D}_2$. 
The function $\mathcal{E}(t,k)$ is time-periodic with period $\tau$  (see \cite{tperiodicI}) and satisfies
$$\mathcal{E}_t + 2ik^2[\sigma_3, \mathcal{E}] = V^b \mathcal{E}.$$
Thus,
$$\mathcal{E}(t,k) = S^b(k) - \int_t^\tau e^{2ik^2(t'-t)\hat{\sigma}_3} (V^b \mathcal{E})(t',k)dt', \qquad 0 \leq t \leq \tau.$$
Since $S^b(k) = I + O(1/k)$, it follows that
$$[\mathcal{E}(t,k)]_1 = \begin{pmatrix} 1 \\ 0 \end{pmatrix} + O(k^{-1}), \qquad k \to \infty, \quad k \in \bar{\mathcal{D}}_2 \cup \bar{\mathcal{D}}_4.$$
Together with the asymptotic formula
\begin{align}\nonumber
[\psi(t,k)]_2 = &\; \Bigg\{\begin{pmatrix} 0 \\ 1 \end{pmatrix} + \frac{1}{k}\begin{pmatrix} - \frac{ig_0^b(t)}{2} \\ i\eta_1(t) \end{pmatrix}
 + \frac{1}{k^2} \begin{pmatrix} \frac{g_1^b(t)}{4} + \frac{g_0^b(t)}{2}\eta_1(t) \\
i\eta_2(t) + \frac{\lambda}{4}|g_0^b(0)|^2
\end{pmatrix} + O\bigg(\frac{1}{k^3}\bigg)\Bigg\}e^{2ik^2t}
	\\ \nonumber
& + \Bigg\{ \frac{1}{k}\begin{pmatrix} \frac{ig_0^b(0)}{2} \\ 0 \end{pmatrix}
 + \frac{1}{k^2} \begin{pmatrix} \frac{g_0^b(0)}{2}\eta_1(t) - \frac{g_1^b(0)}{4} \\ - \frac{\lambda}{4} \overline{g_0^b(t)} g_0^b(0) \end{pmatrix}  + O\bigg(\frac{1}{k^3}\bigg)\Bigg\} e^{-2ik^2t},
 	\\\label{psi2expansion}
& \hspace{7cm} k \to \infty, \quad k \in \C,	
\end{align}
this implies that $m^{(t)}(t, \cdot) \in I + \dot{E}^2(\mathcal{D}_j)$ for $j = 2,4$. Similarly, $m^{(t)}(t, \cdot) \in I + \dot{E}^2(\mathcal{D}_j)$ for $j = 1,3$.
Hence $m^{(t)}(t,\cdot)$ satisfies the $L^2$-RH problem
\begin{align}\label{mtRH}
\begin{cases}
m^{(t)}(t, \cdot) \in I + \dot{E}^2(\C \setminus \Gamma), \\
m^{(t)}_+(t,k) = m^{(t)}_-(t, k)J^{(t)}(t,k)   \quad \text{for a.e.} \ k \in \Gamma,
\end{cases}
\end{align}
where
$$J^{(t)}(t,k) = \begin{pmatrix} 1 & \frac{B^b(k)}{A^b(k)} e^{-4ik^2t} \\
- \frac{\overline{B^b(\bar{k})}}{\overline{A^b(\bar{k})}} e^{4ik^2t} & \frac{1}{A^b(k)\overline{A^b(\bar{k})}} \end{pmatrix}.$$
Moreover, the asymptotic formula (\ref{psi2expansion}) implies that
\begin{align}\label{gjlim}
\begin{cases}
 g_0^b(t) = 2i (m_1)_{12}(t),
	\\
 g_1^b(t) = 4(m_2)_{12}(t) + 2ig_0^b(t) (m_1)_{22}(t),
\end{cases}
\end{align}
where 
$$m_j(t) = \lim_{k \to \infty} k^j m^{(t)}(t,k), \qquad j = 1,2.$$

On the other hand, the relation (\ref{GRquotient}) implies that the function $d(k)$  defined by
$$d(k) = a(k)\overline{A^b(\bar{k})} - b(k) \overline{B^b(\bar{k})}$$
satisfies
$$d(k) = a(k)\overline{A^b(\bar{k})} - \frac{B^b(k)a(k)}{A^b(k)}\overline{B^b(\bar{k})} = \frac{a(k)}{A^b(k)}, \qquad \im k \geq 0.$$
In particular, by (\ref{hApm1}),
$$\bigg(\frac{d(k)}{\overline{A^b(\bar{k})}}\bigg)^{\pm 1}  \in 1 + \dot{E}^2(\mathcal{D}_2) \cap E^\infty(\mathcal{D}_2), \qquad
\bigg(\frac{\overline{d(\bar{k})}}{A^b(k)}\bigg)^{\pm 1}  \in 1 + \dot{E}^2(\mathcal{D}_3) \cap E^\infty(\mathcal{D}_3),$$
which shows that the functions
\begin{align}\nonumber
& G_1(t,k) =  \begin{pmatrix} a(k) & 0 \\ 0 & \frac{1}{a(k)} \end{pmatrix}, 
\qquad G_2(t,k) = \begin{pmatrix} \frac{d(k)}{\overline{A^b(\bar{k})}} & - b(k) e^{-4ik^2t} \\ 0 & \frac{\overline{A^b(\bar{k})}}{d(k)} \end{pmatrix}, 
	\\ \label{Gjdef}
& G_3(t,k) = \begin{pmatrix} \frac{A^b(k)}{\overline{d(\bar{k})}} & 0 \\ - \overline{b(\bar{k})} e^{4ik^2t} & \frac{\overline{d(\bar{k})}}{A^b(k)} \end{pmatrix},
\qquad
G_4(t,k) = \begin{pmatrix} \frac{1}{\overline{a(\bar{k})}} & 0 \\ 0 & \overline{a(\bar{k})} \end{pmatrix},
\end{align}
satisfy $G_j(t,\cdot) \in \dot{E}^2(\mathcal{D}_j) \cap E^\infty(\mathcal{D}_j)$ for $j = 1,\dots, 4$.
Thus a contour deformation argument together with the expression (\ref{Jdef}) for $J$ show that the function $M^{(t)}(t,k)$ defined by
$$M^{(t)}(t,k) = M(0,t,k)G_j(t,k), \qquad k \in \mathcal{D}_j, \quad j = 1,\dots,4,$$
also satisfies (\ref{mtRH}). 
By uniqueness, $m^{(t)} = M^{(t)}$. 
Now
\begin{align}\label{gjlim2}
\begin{cases}
 u(0,t) = 2i (M_1)_{12}(t),
	\\
 u_x(0,t) = 4(M_2)_{12}(t) + 2iu(0,t) (M_1)_{22}(t),
\end{cases}
\end{align}
where 
$$M_j(t) = \lim_{k \to \infty} k^j M^{(t)}(t,k), \qquad j = 1,2.$$
Comparing (\ref{gjlim}) and (\ref{gjlim2}), the claim follows.
\proofendcontinue

\begin{remark}\upshape
The RH problem (\ref{mtRH}) is different from the one used in the analogous situation in \cite{FIS2005}. The formulation in (\ref{mtRH}) has the advantage for us that the matrices $G_1$ and $G_4$ only involve the functions $a(k)$ and $\overline{a(\bar{k})}$, which are analytic in $\mathcal{D}_1$ and $\mathcal{D}_4$, respectively. The matrices $G_2$ and $G_3$ in (\ref{Gjdef}) involve $d(k)$ and $A^b(k)$ as well as their Schwartz conjugates, but since the domains $\mathcal{D}_2$ and $\mathcal{D}_3$ have been defined to be free from branch points and branch cuts, the matrices $G_2$ and $G_3$ still possess the appropriate analyticity properties. 
\end{remark}

\medskip\noindent
{\bf Claim 10.} $u(x,t)$ is given by (\ref{ulimMhat}).

{\it Proof of Claim 10.}
Let 
$$\hat{M} = \begin{cases} 
M \begin{pmatrix} 1 & 0 \\ -h(k) e^{2i(kx + 2k^2t)} & 1 \end{pmatrix}, & k \in \mathcal{D}_1,
	\\
M \begin{pmatrix} 1 & -\overline{h(\bar{k})} e^{-2i(kx + 2k^2t)} \\ 0 & 1 \end{pmatrix}, & k \in \mathcal{D}_4,
	\\
M, &	k \in \mathcal{D}_2 \cup \mathcal{D}_3.
\end{cases}$$
Then the function $\hat{M}$ has no jump discontinuities, but it has singularities at the poles of $h(k)$. Let $\{k_j\}_1^N$ denote the poles of $h(k)$ in $i\R_+$. Then $\hat{M}$ satisfies the residue conditions:
\begin{align*}
& \underset{k = k_j}{\res} [\hat{M}(x,t,k)]_1 = -h_j e^{2i(k_jx + 2k_j^2t)} [\hat{M}(x,t,k_j)]_2,
	\\
& \underset{k = \bar{k}_j}{\res} [\hat{M}(x,t,k)]_2 = -\bar{h}_j e^{-2i(\bar{k}_jx + 2\bar{k}_j^2t)} [\hat{M}(x,t,\bar{k}_j)]_1,
\end{align*}
where $h_j$ denotes the residue of $h(k)$ at $k_j$. This RH problem can be solved using Darboux transformations (see Proposition 2.4 in \cite{FI1996}) and the solution is given by (\ref{Mhatexplicit}).
\proofendcontinue

\medskip\noindent
{\bf Claim 11.}  $u(x,t)$ is periodic in $t$ with period $\tau$ for each $x \geq 0$. In particular, $u(\cdot, t) \in \mathcal{S}([0,\infty))$ for each $t \geq 0$ and $\|u(\cdot, t)\|_{L^1([0,\infty))}$ is bounded for $t \geq 0$.

{\it Proof of Claim 11.}
The fact that all poles of $h(k)$  belong to the set (\ref{hpoleset}) implies that each of the functions $d_j(x,t)$ in (\ref{Bjdef}), and hence also $u(x,t)$, is periodic in $t$ with period $\tau$. 
The function $u(x,t)$ is a solution of (\ref{nls}) on the half-line with initial data in $\mathcal{S}([0,\infty))$ and smooth periodic boundary values. The claim follows from the time-periodicity and functional analytic arguments, 
or from an application of the nonlinear steepest descent method \cite{DZ1993} to the RH problem (\ref{RHm}) which ensures that $u(x,t)$ decays as $x \to \infty$. 
\proofend

\section{Proof of Theorem \ref{expth}}\nequation\label{proof2sec}
Suppose $\{g_0^b(t), g_1^b(t)\} = \{\alpha e^{i\omega t}, ce^{i\omega t}\}$, where $\alpha>0$, $\omega \in \R$, and $c \in \C$. Then (see \cite{Ldefocusing})
$$ S^b(k) = \sqrt{\frac{2\Omega(k) - H(k)}{2\Omega(k)}} \begin{pmatrix} 1 & \frac{iH(k)}{2\alpha k - i\bar{c}} \\
-\frac{iH(k)}{2\alpha k + ic} & 1 \end{pmatrix},$$	
where
\begin{align*}
 \Omega(k) = \sqrt{4k^4 + 2\omega k^2 + 4 \alpha \im(c) k + \bigg(\frac{\omega}{2} +  \alpha^2\bigg)^2 - |c|^2},
\end{align*}
and 
$$H(k) = \Omega(k) - 2k^2 - \alpha^2 - \frac{\omega}{2}.$$
In particular, 
$$Q^b(k) = \frac{B^b(k)}{A^b(k)} = \frac{iH(k)}{2\alpha k - i\bar{c}}.$$
The set of branch points consists of the four zeros of $\Omega^2(k)$  and the points $\{\frac{i\bar{c}}{2\alpha}, -\frac{ic}{2\alpha}\}$. Indeed, the identity
\begin{align}\label{2OmegaHH}
  (H - 2\Omega)H = (2\alpha k - i\bar{c})(2\alpha k + ic)
\end{align}  
implies that the zeros of $2\Omega - H$ and $H$ are included in $\{\frac{i\bar{c}}{2\alpha}, -\frac{ic}{2\alpha}\}$.

The result of \cite{Ldefocusing} implies that if the pair $\{\alpha e^{i\omega t}, ce^{i\omega t}\}$ is eventually admissible for the defocusing NLS equation, then the parameter triple $(\alpha, \omega, c)$ belongs to one of the following families:
\begin{subequations}\label{admissiblesets}
\begin{align}\label{admissibleA}
& \bigg\{\bigg(\alpha, \omega, c = \pm \sqrt{\frac{(\omega + 3\alpha^2)^3}{27\alpha^2}} + \frac{i|\omega|^{3/2}}{3\sqrt{3} \alpha}\bigg) \; \bigg| \;  \alpha > 0, \; -3\alpha^2 \leq \omega < 0\bigg\},
	\\ \nonumber
& \bigg\{\bigg(\alpha  = -\frac{4K^3 + \omega K}{c_2}, \omega, c = \pm \sqrt{ \bigg(\alpha^2+ \frac{\omega}{2}\bigg)^2 -c_2^2 - 2K^2(6K^2 + \omega)} + ic_2\bigg)  
	\\ \label{admissibleB}
&\hspace{3.5cm} \bigg| \;  -12 K^2 < \omega < -4K^2, \; 0 < c_2 \leq -\frac{4K^2 + \omega}{2}, \; K > 0\bigg\}, 
	\\ \label{admissibleC}
& \big\{(\alpha,\omega, c = i\alpha\sqrt{-2\alpha^2 - \omega}) \; \big| \; \alpha > 0, \; \omega < -3\alpha^2\big\},
	\\ \label{admissibleD}
& \big\{(\alpha,\omega,c = \pm \alpha\sqrt{\omega + \alpha^2}) \; \big| \; \omega + \alpha^2 \geq 0, \; \alpha > 0\big\},
	\\ \nonumber
& \bigg\{\bigg(\alpha = -\frac{4K^3 + \omega K}{c_2}, \omega, c = \pm \sqrt{\bigg(\alpha^2+ \frac{\omega}{2}\bigg)^2 -c_2^2 - 2K^2(6K^2 + \omega)} + ic_2\bigg)
	\\ \label{admissibleE}
& \hspace{3.5cm} \bigg| \; -4K^2 < \omega \leq -3K^2, \; -\frac{4K^2 + \omega}{2} \leq c_2 < 0, \; K > 0\bigg\}.
\end{align}
\end{subequations}

We will consider each of the families in (\ref{admissiblesets}) in turn. Theorem \ref{mainth} will show that the only eventually admissible pairs arise from the subfamily of (\ref{admissibleD}) given in (\ref{reallyadmissibleD}).

\subsection{Family (\ref{admissibleA})}
\begin{figure}
\begin{center}
\begin{overpic}[width=.45\textwidth]{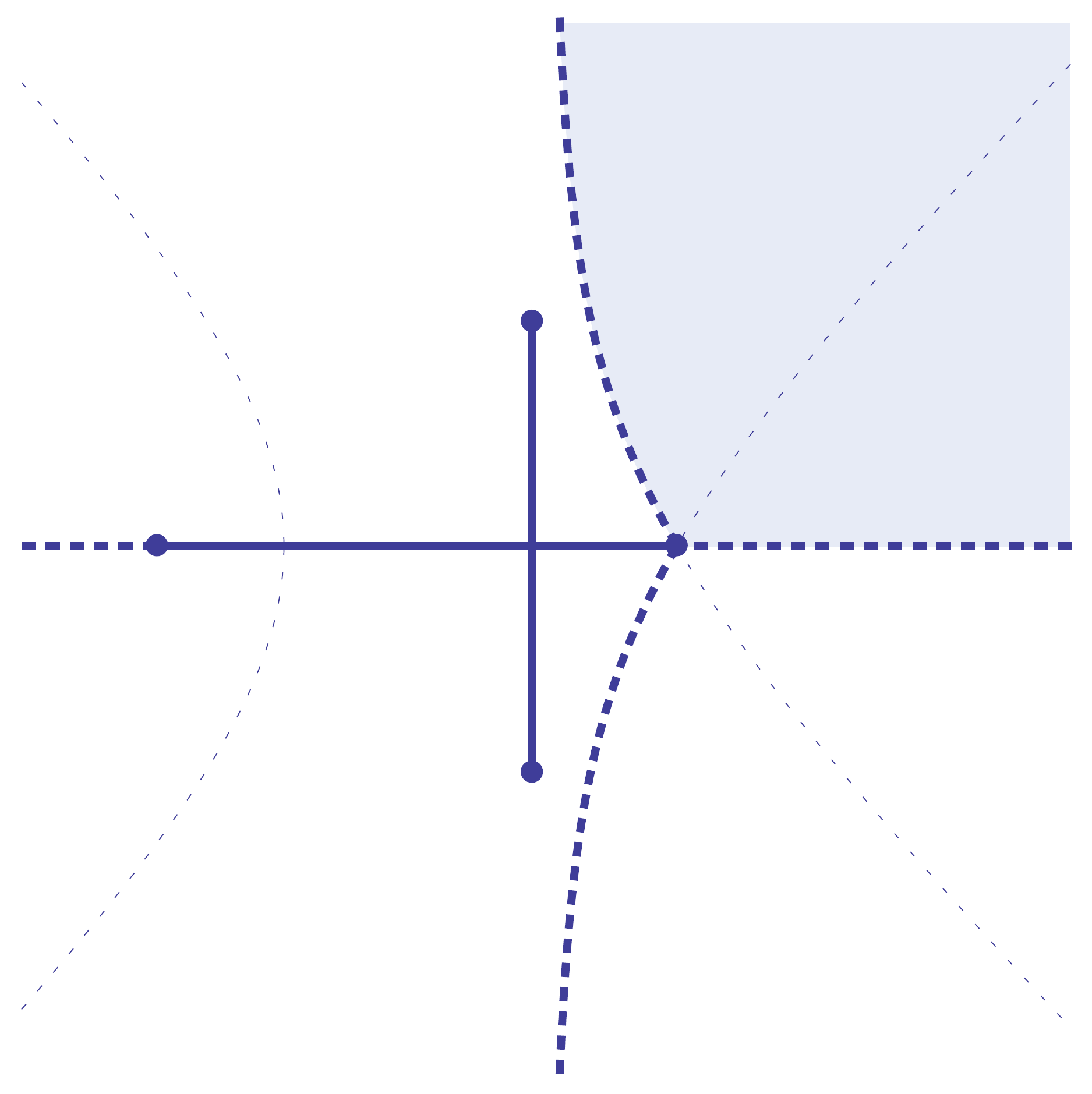}
      \put(65,75){$D_1$}
      \put(28,69){$\frac{c_2 + i|c_1|}{2\alpha}$}
      \put(28,28){$\frac{c_2 - i|c_1|}{2\alpha}$}
      \put(61,43){$K$}
      \put(6,43){$-3K$}
\end{overpic}
     \begin{figuretext}\label{cutsA.pdf}
       The  branch cuts for the triples in (\ref{admissibleA}). Branch points are indicated by dots, branch cuts by solid lines, $\im \Omega(k) = 0$ on the thick striped curves, and $\re \Omega(k) = 0$ on the thin striped curves.
      \end{figuretext}
     \end{center}
\end{figure}

For the family (\ref{admissibleA}), we have
$$\Omega(k) = 2(k-K)\sqrt{(k-K)(k+3K)}$$
where $K = \sqrt{\frac{|\omega|}{12}}$, see Figure \ref{cutsA.pdf}. Since $\Omega'(k)$, and hence also $(Q^b)'(k)$, is singular at $-3K$, Theorem \ref{mainth} implies that the pairs in this family are not eventually admissible.

\subsection{Family (\ref{admissibleB})}
\begin{figure}
\begin{center}
\begin{overpic}[width=.45\textwidth]{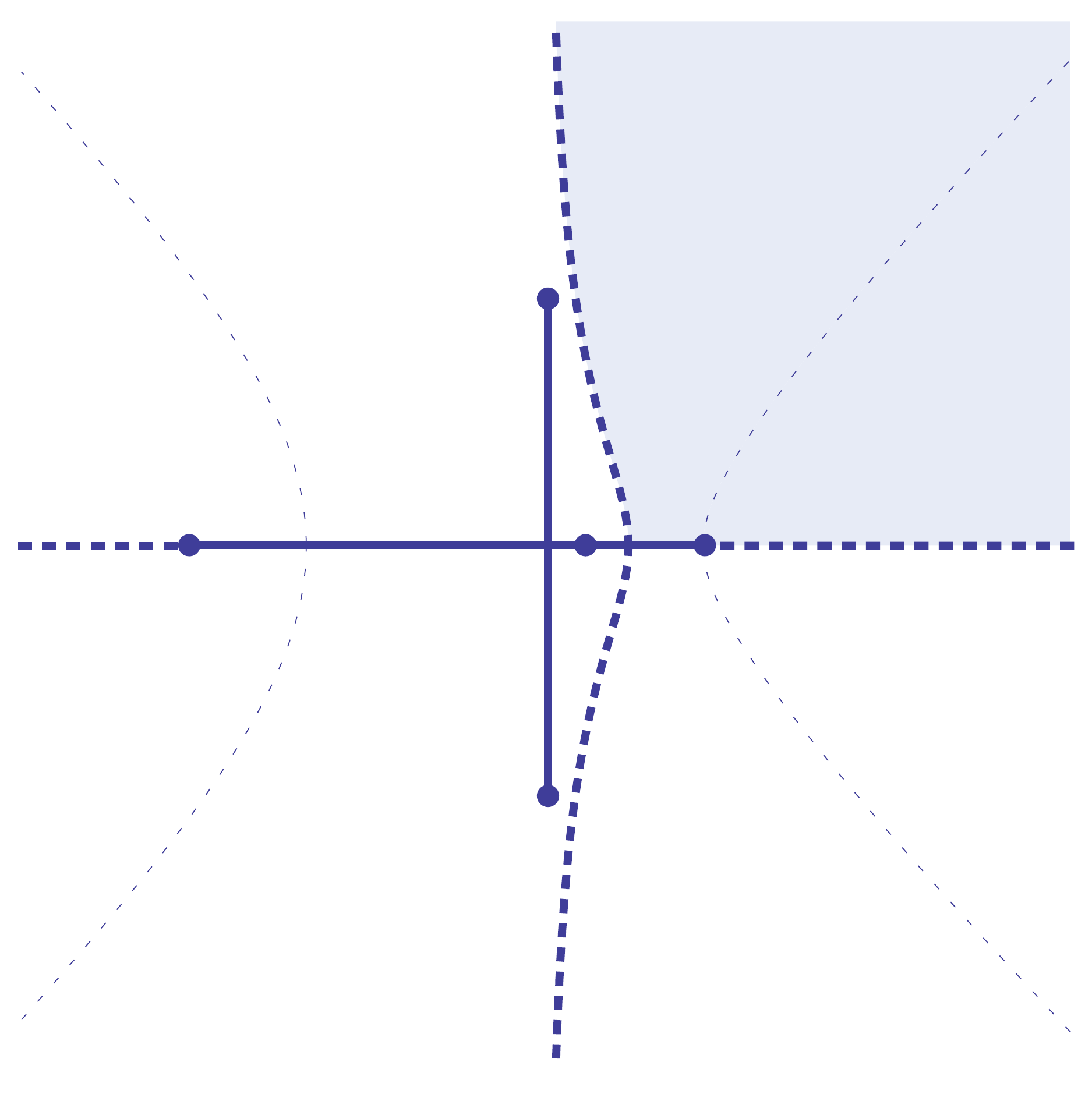}
      \put(65,75){$D_1$}
      \put(30,71){$\frac{c_2 + i|c_1|}{2\alpha}$}
      \put(30,26){$\frac{c_2 - i|c_1|}{2\alpha}$}
      \put(15,43){$K_1$}
      \put(51,43){$K_2$}
      \put(61,43){$K$}
\end{overpic}
     \begin{figuretext}\label{cutsB.pdf}
       The  branch cuts for the triples in (\ref{admissibleB}). 
      \end{figuretext}
     \end{center}
\end{figure}

For the family (\ref{admissibleB}), we have
$$\Omega(k) = 2(k-K)\sqrt{(k-K_1)(k-K_2)},$$
where (see Figure \ref{cutsB.pdf})
$$K_1 = -K-\sqrt{-2 K^2-\frac{\omega }{2}}, \qquad K_2 = -K + \sqrt{-2 K^2-\frac{\omega }{2}}.$$
Since $\Omega'(k)$, and hence also $(Q^b)'(k)$, is singular at $K_1$ and $K_2$, Theorem \ref{mainth} implies that the pairs in this family are not eventually admissible.

\subsection{Family (\ref{admissibleC})}
\begin{figure}
\begin{center}
\begin{overpic}[width=.45\textwidth]{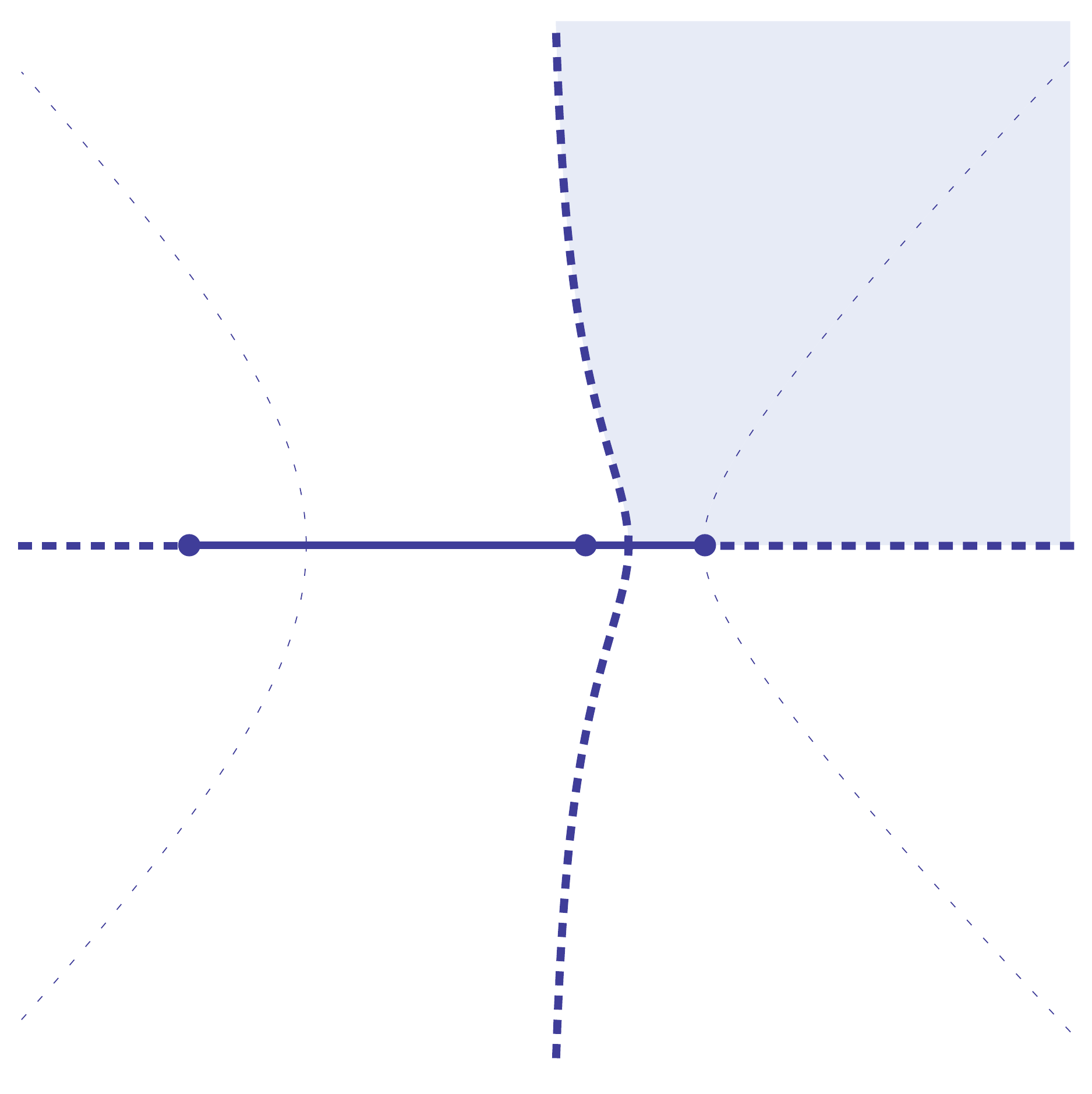}
      \put(65,80){$D_1$}
      \put(62,43){$K$}
   \put(14,43){$K_1$}
   \put(47,43){$K_2$}
\end{overpic}
     \begin{figuretext}\label{cutsC.pdf}
        The  branch cuts for the triples in (\ref{admissibleC}). Since $c_1 = 0$, the branch cut $(-\frac{ic}{2\alpha}, \frac{i\bar{c}}{2\alpha})$ is absent. 
     \end{figuretext}
     \end{center}
\end{figure}

For the family (\ref{admissibleC}), we have
$$\Omega(k) = 2(k-K)\sqrt{(k-K_1)(k - K_2)},$$
where (see Figure \ref{cutsC.pdf})
$$K_1 = -K-\alpha, \qquad K_2 = -K + \alpha.$$ 
Since $\Omega'(k)$, and hence also $(Q^b)'(k)$, is singular at $K_1$ and $K_2$, Theorem \ref{mainth} implies that the pairs in this family are not eventually admissible.

\subsection{Family (\ref{admissibleD})}
For the family (\ref{admissibleD}), we have
$$\Omega(k) = 2k^2 + \frac{\omega}{2}, \qquad H(k) = - \alpha^2.$$
Hence
$$Q^b(k) = \begin{cases} \frac{-i\alpha}{2 k - i\sqrt{\omega + \alpha^2}} &\text{if} \quad c = \alpha \sqrt{\omega + \alpha^2}, \\
 \frac{-i\alpha}{2 k + i\sqrt{\omega + \alpha^2}} &\text{if} \quad c = -\alpha \sqrt{\omega + \alpha^2}.
\end{cases}$$
The pairs with $c = \alpha \sqrt{\omega + \alpha^2}$ are inadmissible because in this case $Q^b(k)$  has a pole in the upper half plane.
The pairs with $c = - \alpha \sqrt{\omega + \alpha^2}$ and $-\alpha^2 \leq \omega \leq 0$ are inadmissible because in this case
$$\sup_{k\in\R} |Q^b(k)|^2 = |Q^b(0)|^2 =  \frac{\alpha^2}{\omega + \alpha^2} \geq 1,$$ 
contradicting property (A3) of Theorem \ref{mainth}.

\begin{figure}
\begin{center}
\begin{overpic}[width=.45\textwidth]{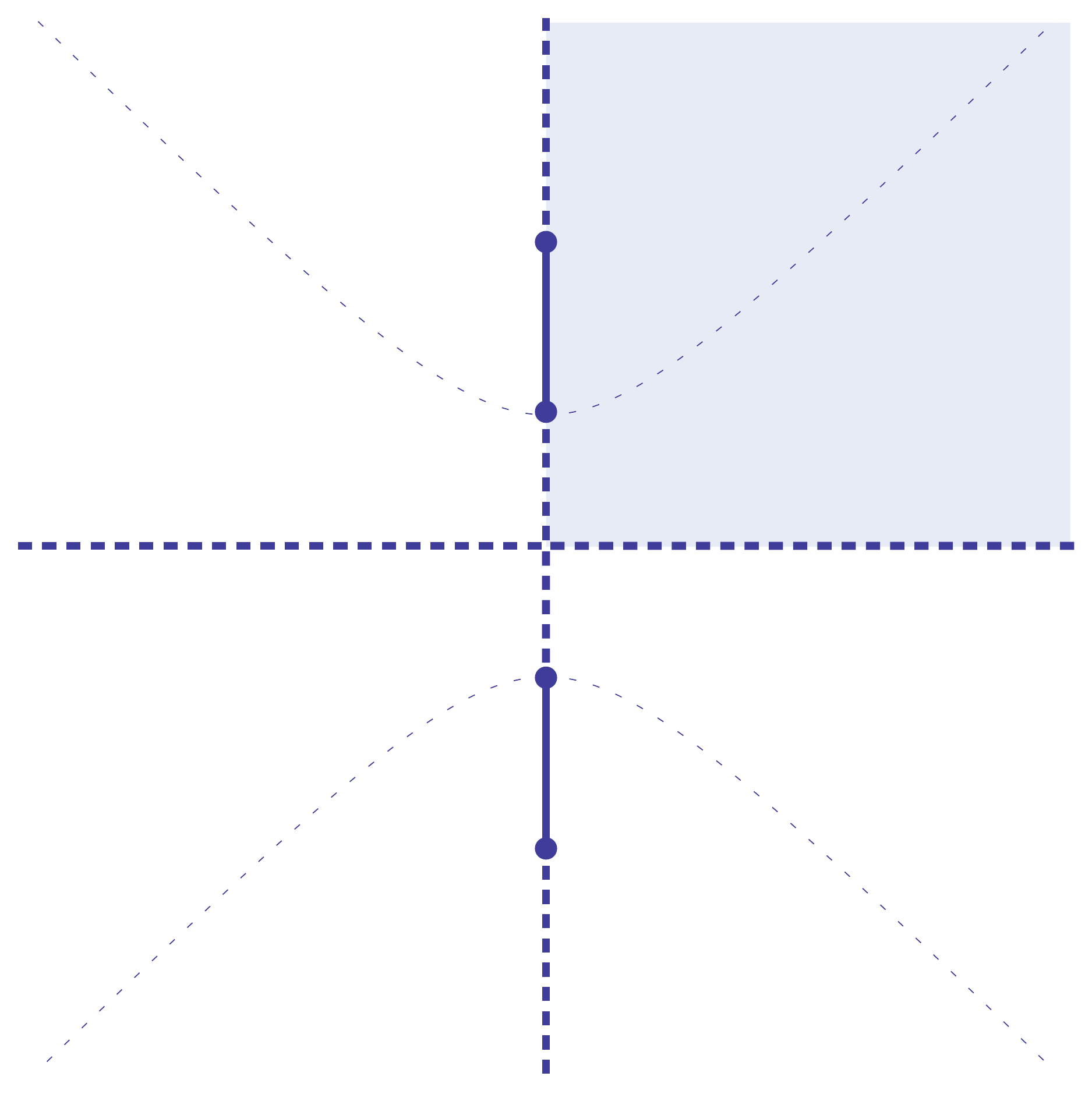}
      \put(75,65){$D_1$}
      \put(37,60){$\frac{i\sqrt{\omega}}{2}$}
      \put(32,37){$-\frac{i\sqrt{\omega}}{2}$}
      \put(52,78){$\frac{i\sqrt{\omega + \alpha^2}}{2}$}
      \put(52,20){$-\frac{i\sqrt{\omega + \alpha^2}}{2}$}
\end{overpic}
     \begin{figuretext}\label{cutsDa.pdf}
         The  branch cuts for the triples in  (\ref{admissibleD}) in the case of $\omega > 0$.    \end{figuretext}
     \end{center}
\end{figure}

\begin{figure}
\begin{center}
\begin{overpic}[width=.45\textwidth]{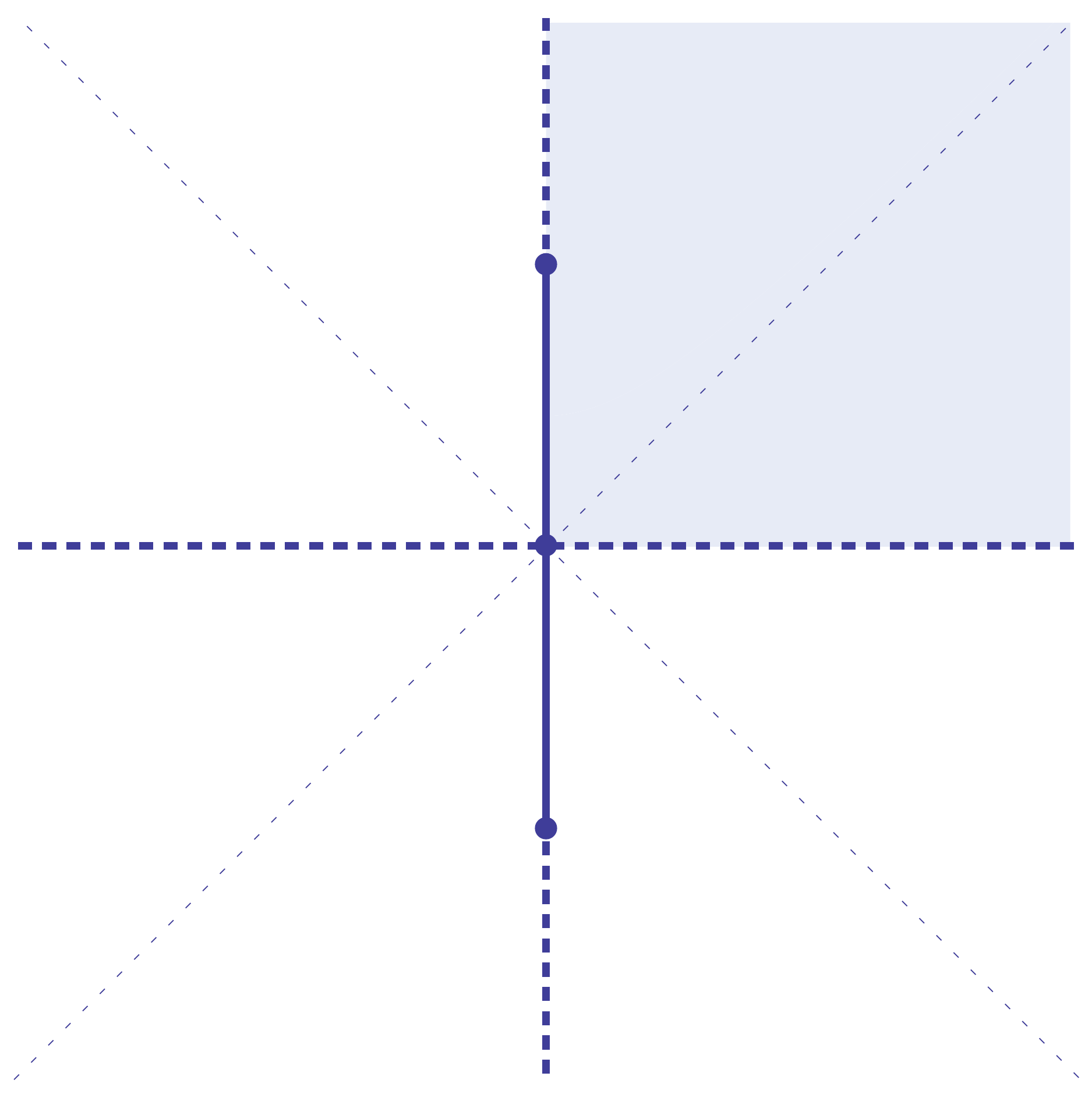}
      \put(75,65){$D_1$}
      \put(52,74){$\frac{i\alpha}{2}$}
      \put(52,23){$-\frac{i\alpha}{2}$}
\end{overpic}
     \begin{figuretext}\label{cutsDb.pdf}
         The  branch cuts for the triples in  (\ref{admissibleD}) in the case of $\omega = 0$.
     \end{figuretext}
     \end{center}
\end{figure}

\begin{figure}
\begin{center}
\begin{overpic}[width=.45\textwidth]{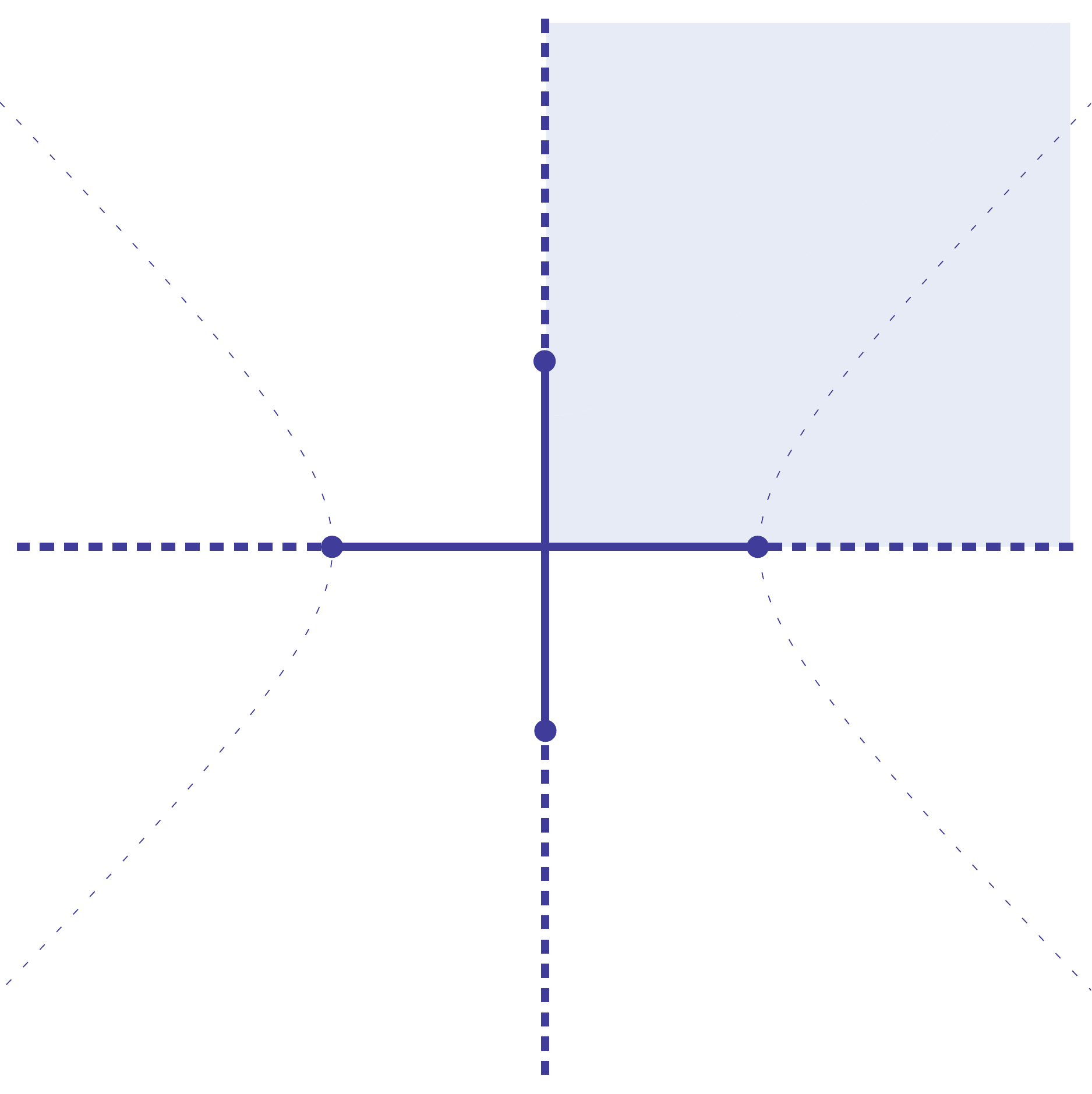}
     \put(75,80){$D_1$}
      \put(66,40){$\frac{\sqrt{|\omega|}}{2}$}
      \put(17,40){$-\frac{\sqrt{|\omega|}}{2}$}
      \put(52,67){$\frac{i\sqrt{\omega + \alpha^2}}{2}$}
      \put(51,29){$-\frac{i\sqrt{\omega + \alpha^2}}{2}$}
\end{overpic}
     \begin{figuretext}\label{cutsDc.pdf}
       The  branch cuts for the triples in  (\ref{admissibleD}) in the case of $-\alpha^2 \leq \omega < 0$.
      \end{figuretext}
     \end{center}
\end{figure}

Thus suppose $c = -\alpha \sqrt{\omega + \alpha^2}$ and $\omega > 0$. In this case, $Q^b(k)$ satisfies the three conditions (A1)-(A3) of Theorem \ref{mainth}.
Moreover, 
$$A^b(k) = \sqrt{\frac{(k+k_2)(k - k_2)}{(k+k_1)(k-k_1)}}, \qquad B^b(k) = \frac{\alpha}{2i} \sqrt{\frac{k - k_2}{(k+k_1)(k-k_1)(k+k_2)}}.$$
where
$$k_1 = \frac{i\sqrt{\omega}}{2}, \qquad k_2 = \frac{i\sqrt{\omega + \alpha^2}}{2}.$$
Hence $A^b(k)^2$ and $P^b(k) = -\frac{\alpha}{2i}\frac{k + k_2}{(k+k_1)(k-k_1)}$ are rational functions with simple poles. 
This leads to
$$a(k) = \frac{k + k_2}{k + k_1}, \qquad b(k) = \frac{\alpha}{2i(k + k_1)}, \qquad d(k) = \sqrt{\frac{(k-k_1)(k + k_2)}{(k+k_1)(k-k_2)}},$$
and 
$$h(k) = \frac{\alpha(k+k_1)}{2i(k-k_1)(k+k_2)}.$$
Following Theorem \ref{mainth}, we define $\hat{M}(x,t,k)$ by
$$\hat{M}(x,t,k) = (kI + B_1(x,t))\begin{pmatrix} \frac{1}{k-k_1} & 0 \\ 0 & \frac{1}{k-\bar{k}_1} \end{pmatrix},$$
where $k_1 = \frac{i\sqrt{\omega}}{2}$ is the only pole of $P^b$ in the upper half plane and the matrix $B_1(x,t)$ is determined from 
\begin{align}\label{B1def}
  (k_1I + B_1(x,t)) \begin{pmatrix} 1 \\ -d_1(x,t) \end{pmatrix} = 0, 
\qquad   (\bar{k}_1I + B_1(x,t)) \begin{pmatrix} - \overline{d_1(x,t)} \\ 1 \end{pmatrix} = 0,
\end{align}
with 
\begin{align*}
& d_1(x,t) = -h_1 \frac{1}{k_1 - \bar{k}_1} e^{2i(k_1x + 2k_1^2t)},
	\\
& h_1 = \underset{k = k_1}{\res} h(k) = -\frac{i\alpha\sqrt{\omega}}{\sqrt{\omega} + \sqrt{\alpha^2 + \omega}}.
\end{align*}
The general solution of (\ref{B1def}) is given by
$$(B_1(x,t))_{12} = -\frac{4 \bar{h}_1 k_1^2 e^{-2 i k_1 (2 k_1 t-x)}}{4 k_1^2+h_1 \bar{h}_1 e^{4 i k_1 x}}, \qquad (B_1(x,t))_{22} = \frac{k_1 \left(4 k_1^2-h_1 \bar{h}_1 e^{4 i k_1 x}\right)}{4 k_1^2+h_1 \bar{h}_1 e^{4 i k_1 x}},$$
showing that the RH problem for $\hat{M}$ is solvable except when
$$4 k_1^2+h_1 \bar{h}_1 e^{4 i k_1 x} = 0, \quad \text{i.e.} \quad x =  -\frac{\log \left(\frac{2 \sqrt{\omega } \sqrt{\alpha ^2+\omega }+\alpha ^2+2 \omega}{\alpha ^2}\right)}{2 \sqrt{\omega }}.$$
In particular, the RH problem for $\hat{M}$ is solvable for all $(x,t) \in [0,\infty) \times [0,\infty)$. 
Theorem \ref{mainth} therefore shows that the pairs with $c = -\alpha \sqrt{\omega + \alpha^2}$ and $\omega > 0$ are admissible and that the unique solution $u(x,t)$ which satisfies (\ref{g0g0Bg1g1B}) is given by
\begin{align}\label{ulimhatM}
u(x,t) = 2i \lim_{k \to \infty} (k\hat{M}(x,t,k))_{12},
\end{align}
where
\begin{align*}
& \hat{M}_{12}(x,t) = -\frac{2 \alpha  \sqrt{\omega } \left(\sqrt{\alpha ^2+\omega }+\sqrt{\omega }\right) e^{x
   \sqrt{\omega }+i t \omega }}{\left(\sqrt{\omega }-2 i k\right) \big(\alpha ^2
   \big(e^{2 x \sqrt{\omega }}-1\big)+2 \big(\sqrt{\omega  (\alpha ^2+\omega)}+\omega \big) e^{2 x \sqrt{\omega }}\big)},
  	\\
& \hat{M}_{22}(x,t) = \frac{1}{k+\frac{i
   \sqrt{\omega }}{2}}
   \left(k+\frac{i \sqrt{\omega } \big(\alpha ^2 \big(e^{2 x \sqrt{\omega }}+1\big)+2
   \big(\sqrt{\omega  \big(\alpha ^2+\omega \big)}+\omega \big) e^{2 x \sqrt{\omega
   }}\big)}{2 \alpha ^2 \big(e^{2 x \sqrt{\omega }}-1\big)+4 \big(\sqrt{\omega 
   \big(\alpha ^2+\omega \big)}+\omega \big) e^{2 x \sqrt{\omega }}}\right).
\end{align*}	
Computing the limit in (\ref{ulimhatM}), we find (\ref{uexplicit}).

It is of course easy to verify directly that the solution in (\ref{uexplicit}) satisfies (\ref{nls}) as well as the correct boundary conditions:
$$u(0,t) = \alpha e^{i\omega t}, \qquad u_x(0,t) = -\alpha \sqrt{\omega + \alpha^2} e^{i\omega t}.$$

\subsection{Family (\ref{admissibleE})}
\begin{figure}
\begin{center}
\begin{overpic}[width=.45\textwidth]{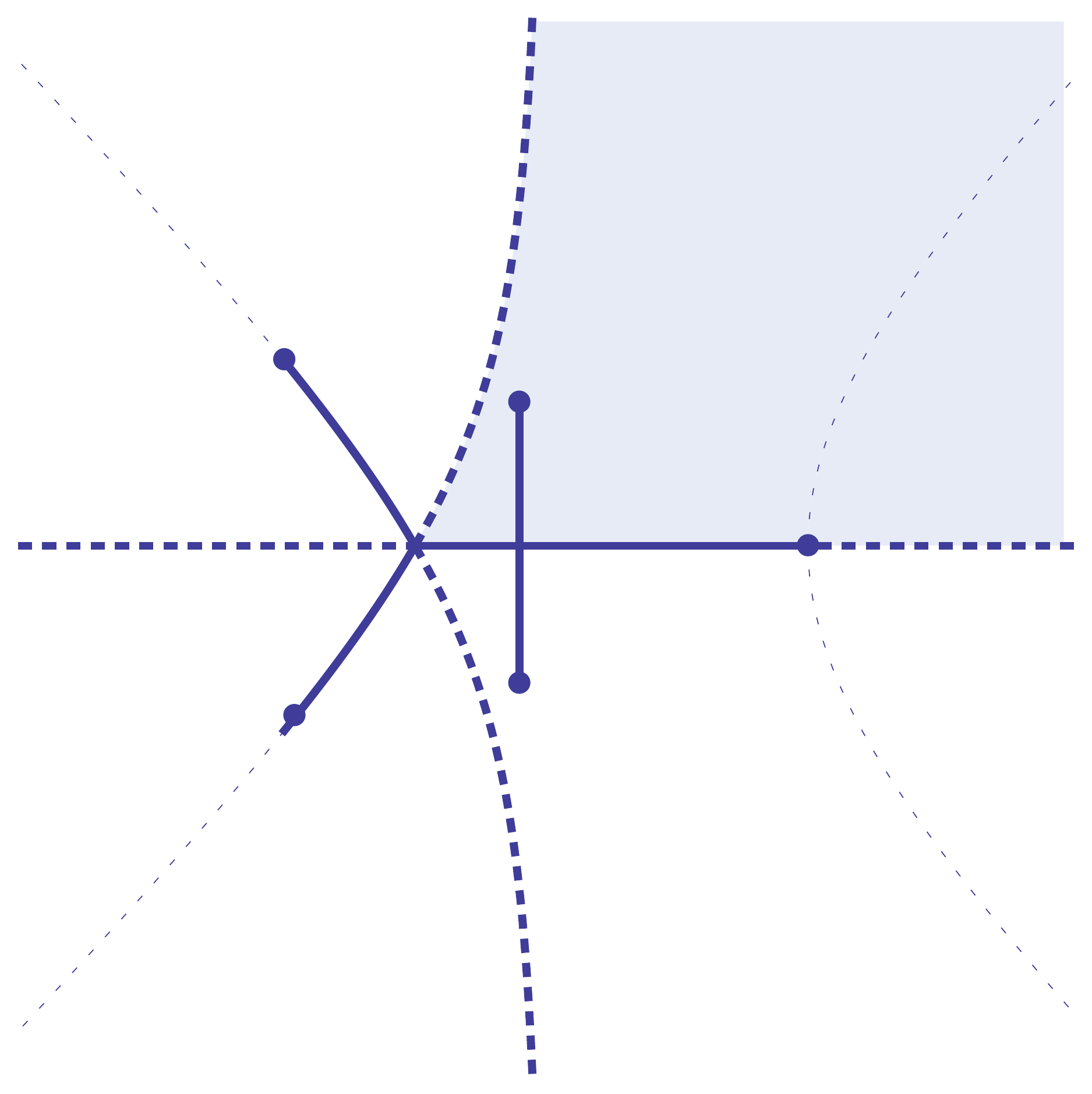}
     \put(65,75){$D_1$}
     \put(75,44){$K$}
      \put(16,63){$K_1$}
     \put(16,35){$\bar{K}_1$}
      \put(50,63){$\frac{c_2 + i|c_1|}{2\alpha}$}
      \put(50,36){$\frac{c_2 - i|c_1|}{2\alpha}$}
\end{overpic}
     \begin{figuretext}\label{cutsEa.pdf}
         The  branch cuts for the triples in  (\ref{admissibleE}) in the case of $\omega = -3K^2$.
     \end{figuretext}
  \end{center}
\end{figure}

\begin{figure}
\begin{center}
\begin{overpic}[width=.45\textwidth]{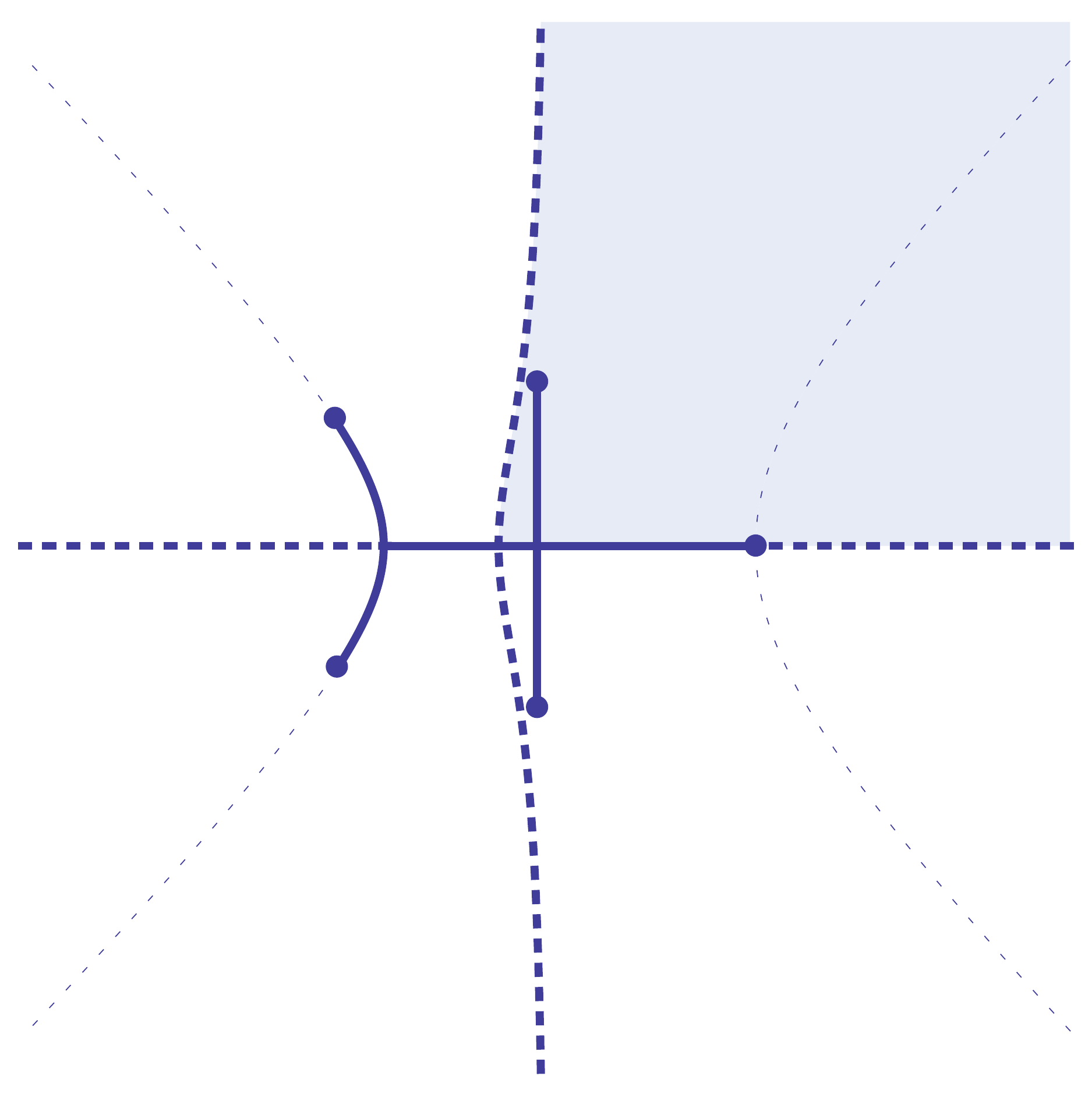}
     \put(65,75){$D_1$}
     \put(70,44){$K$}
      \put(21,60){$K_1$}
     \put(21,38){$\bar{K}_1$}
      \put(51,64){$\frac{c_2 + i|c_1|}{2\alpha}$}
      \put(51,33){$\frac{c_2 - i|c_1|}{2\alpha}$}
\end{overpic}
     \begin{figuretext}\label{cutsEb.pdf}
         The  branch cuts for the triples in  (\ref{admissibleE}) in the case of $-4K^2 < \omega < -3K^2$. 
     \end{figuretext}
  \end{center}
\end{figure}

For the family (\ref{admissibleE}), we have
$$\Omega(k) = 2(k-K)\sqrt{(k-K_1)(k-\bar{K}_1)}$$
where $K_1 = (-1 + \frac{i}{\sqrt{2}})K$ if $\omega = -3K^2$ (see Figure \ref{cutsEa.pdf}) and $K_1 = -K + i\sqrt{2K^2 + \frac{\omega}{2}}$ if  $-4K^2 < \omega < -3K^2$ (see Figure \ref{cutsEb.pdf}).
Since $\Omega'(k)$, and hence also $(Q^b)'(k)$, is singular at $K_1$, Theorem \ref{mainth} implies that the pairs in this family are not eventually admissible. 
This completes the proof of Theorem \ref{expth}.
\proofend

\section{New exact solutions}\nequation
In Theorem \ref{expth} we derived an explicit solution of (\ref{nls}) in the quarter plane with periodic Dirichlet and Neumann values given by single exponentials.
By applying the procedure of Theorem \ref{mainth} with other choices of the simple poles $\{k_j\}$ and the associated residues $\{h_j\}$, we can generate other exact solutions of the defocusing NLS on the half-line with $t$-periodic boundary values. These solutions bear certain similarities with stationary solitons. They are smooth on the half-line, but singular on the line. 

\subsection{Example}
Let 
$$k_1 = \frac{i}{2}, \qquad k_2 = i, \qquad c_1 = 1, \qquad c_2 = 1+i,$$
Then the procedure of Theorem \ref{expth} yields the following exact solution of (\ref{nls}) in the quarter plane:
$$u(x,t) = \frac{u_1(x,t)}{u_2(x,t)},$$
where
\begin{align*}
 u_1(x,t) = &\; 72 e^{2 (x+8 i t)} \left((72+72 i) e^{6 (x+2 i t)}-(2+2 i) e^{2 (x+6 i t)}+72 i
   e^{8 x}-i\right),
   	\\
 u_2(x,t) = & -36 e^{4 x} \left(18 e^{4 (x+3 i t)}+(8-8 i) e^{2 x+24 i t}+9 e^{12
   i t}+(8+8 i) e^{2 x}\right)
   	\\
&   +2592 e^{12 (x+i t)}+e^{12 i t}.
\end{align*}   
Note that $u(x,t)$ is $t$-periodic for each $x \geq 0$. Moreover, $u(x,t)$ is singular at 
$$x \approx -1.47 \text{ and } x \approx -0.0908.$$

\bigskip
\noindent
{\bf Acknowledgement} {\it The author acknowledges support from the EPSRC, UK.}

\bibliographystyle{plain}
\bibliography{is}

\begin{thebibliography}{99}
\small

\bibitem{BDT1988}
R. Beals, P. Deift, and C. Tomei, {\it Direct and inverse scattering on the line}, Mathematical Surveys and Monographs 28, American Mathematical Society, Providence, RI, 1988.

\bibitem{BFS2003}
A. Boutet de Monvel, A. S. Fokas, and D. Shepelsky, The analysis of the
global relation for the nonlinear Schr\"odinger equation on the half-line,
{\it Lett. Math. Phys.} {\bf 65} (2003), 199--212.

\bibitem{BIK2009}
A. Boutet de Monvel, A. Its, and V. Kotlyarov, Long-time asymptotics for the focusing NLS equation with time-periodic boundary condition on the half-line, {\it Comm. Math. Phys.} {\bf 290} (2009), 479--522.

\bibitem{BK2007}
A. Boutet de Monvel and V. Kotlyarov, The focusing nonlinear Schr\"odinger equation on the quarter plane with time-periodic boundary condition: a Riemann-Hilbert approach, {\it J. Inst. Math. Jussieu} {\bf 6} (2007), 579--611. 

\bibitem{BKS2009}
A. Boutet de Monvel, V. Kotlyarov, and D. Shepelsky, Decaying long-time asymptotics for the focusing NLS equation with periodic boundary condition, {\it Int. Math. Res. Not. IMRN} {\bf 2009}, 547--577.

\bibitem{D1999}
P. Deift, {\it Orthogonal polynomials and random matrices: a Riemann-Hilbert approach}, Courant Lecture Notes in Mathematics, 3, New York University, Courant Institute of Mathematical Sciences, New York; American Mathematical Society, Providence, RI, 1999.

\bibitem{DZ1993}
P. Deift and X. Zhou, A steepest descent method for oscillatory Riemann-Hilbert problems. Asymptotics for the MKdV equation, 
{\it Ann. of Math.} {\bf 137} (1993), 295--368.

\bibitem{DZ2002a}
P. Deift and X. Zhou, Perturbation theory for infinite-dimensional integrable systems on the line. A case study, {\it Acta Math.} {\bf 188} (2002), 163--262.

\bibitem{DZ2002b}
P. Deift and X. Zhou, A priori $L^p$-estimates for solutions of Riemann-Hilbert problems, {\it Int. Math. Res. Not.} {\bf 2002}, 2121--2154.

\bibitem{F1997}
A. S. Fokas, A unified transform method for solving linear and certain nonlinear PDEs, 
{\it Proc. Roy. Soc. Lond.} A {\bf 453} (1997), 1411--1443.

\bibitem{F2002}
A. S. Fokas, Integrable nonlinear evolution equations on the half-line, 
{\it Comm. Math. Phys.} {\bf 230} (2002), 1--39.

\bibitem{F2005}
A. S. Fokas, A generalised Dirichlet to Neumann map for certain nonlinear evolution PDEs, 
{\it Comm. Pure Appl. Math.} {\bf LVIII} (2005), 639--670.

\bibitem{FI1996}
A. S. Fokas and A. R. Its, The linearization of the initial-boundary value problem of the nonlinear Schr\"odinger equation,
{\it SIAM J. Math. Anal.} {\bf 27} (1996), 738--764. 

\bibitem{FIKN2006}
A. S. Fokas, A. R. Its, A. A. Kapaev, V. Y. Novokshenov, {\it Painlev\'e transcendents. The Riemann-Hilbert approach.} Mathematical Surveys and Monographs, 128. American Mathematical Society, Providence, RI, 2006.

\bibitem{FIS2005}
A. S. Fokas, A. R. Its, and L.-Y. Sung, The nonlinear Schr\"odinger equation on the half-line, 
{\it Nonlinearity} {\bf 18} (2005), 1771--1822.

\bibitem{trilogy1}
A. S. Fokas and J. Lenells, The unified method: I. Non-linearizable problems on the half-line, {\it J. Phys. A: Math. Theor.} {\bf 45}, 195201.

\bibitem{LCarleson}
J. Lenells, Matrix Riemann-Hilbert problems with jumps across Carleson contours, preprint, arXiv:1401.2506.
 
 \bibitem{Ldefocusing}
J. Lenells, Admissible boundary values for the defocusing nonlinear Schr\"odinger equation with asymptotically $t$-periodic data, preprint, arXiv:1407.5046.

\bibitem{tperiodicI}
J. Lenells and A. S. Fokas, The nonlinear Schr\"odinger equation with $t$-periodic data: I. Exact results, preprint, arXiv:1412.0304.

\bibitem{tperiodicII}
J. Lenells and A. S. Fokas, The nonlinear Schr\"odinger equation with $t$-periodic data: II. Perturbative results, preprint, arXiv:1412.0306.
 

\bibitem{L1996}
B. Ya. Levin, {\it Lectures on entire functions}, Translations of Mathematical Monographs 150, American Mathematical Society, Providence, RI, 1996.

\bibitem{Z1989}
X. Zhou, The Riemann-Hilbert problem and inverse scattering, {\it SIAM J. Math. Anal.} {\bf 20} (1989), 966--986. 

\end{thebibliography}

\end{document}